\documentclass[11pt]{article}
\usepackage[a4paper, total={6in, 10in}]{geometry}
\usepackage[utf8]{inputenc}
\usepackage{amsthm,mathtools,amsmath}
\usepackage{enumitem}
\usepackage{hyperref}
\usepackage[charter]{mathdesign}
\usepackage{color}
\usepackage{titlesec}
\usepackage{subcaption}

 \usepackage{ulem}

\titleformat{\paragraph}
{\normalfont\normalsize\bfseries}{\theparagraph}{1em}{}
\titlespacing*{\paragraph}
{0pt}{3.25ex plus 1ex minus .2ex}{1.5ex plus .2ex}

\DeclareFontFamily{U}{BOONDOX-calo}{\skewchar\font=45 }
\DeclareFontShape{U}{BOONDOX-calo}{m}{n}{
  <-> s*[1.05] BOONDOX-r-calo}{}
\DeclareFontShape{U}{BOONDOX-calo}{b}{n}{
  <-> s*[1.05] BOONDOX-b-calo}{}
\DeclareMathAlphabet{\mathcalboondox}{U}{BOONDOX-calo}{m}{n}
\SetMathAlphabet{\mathcalboondox}{bold}{U}{BOONDOX-calo}{b}{n}
\DeclareMathAlphabet{\mathbcalboondox}{U}{BOONDOX-calo}{b}{n}
\newcommand{\mcb}[1]{{\mathcalboondox #1}}

\makeatletter \@addtoreset{equation}{section}
\@addtoreset{figure}{section}
 \makeatother

\newtheorem{theorem}{Theorem}[section]
\newtheorem{corollary}{Corollary}[theorem]

\newtheorem{definition}{Definition}[section]
\newtheorem{remark}{Remark}[]

\usepackage{tikz}
\usetikzlibrary{shapes.misc}
\usetikzlibrary{arrows.meta}
\tikzset{cross/.style={cross out, draw=black, minimum size=2*(#1-\pgflinewidth), inner sep=0pt, outer sep=0pt},
cross/.default={0.2cm}}
\definecolor{blue-violet}{rgb}{0.54, 0.17, 0.89}
\usetikzlibrary{decorations.pathreplacing}

\newcommand{\lsim}{\lesssim}

\newenvironment{acknowledgements}{%
  \begin{abstract}
}{%
  \end{abstract}
}
\title{On the correlations  of some  microscopic random systems} 
\author{ P. Gon\c calves\thanks{Instituto Superior T\'ecnico, Department of Mathematics, Av. Rovisco Pais 1, 1049-001, Lisbon. E-mail: {\tt pgoncalves@tecnico.ulisboa.pt}}\and B. Salvador\thanks{Instituto Superior T\'ecnico, Department of Mathematics, Av. Rovisco Pais 1, 1049-001, Lisbon. E-mail: {\tt beatriz.salvador@tecnico.ulisboa.pt}}}

\date{ }

\begin{document}

\maketitle

\begin{abstract}
We investigate the two-points correlation function for several boundary-driven interacting particle systems. Our goal is to show that the time evolution of that correlation function is solution to a partial differential equation that can be written in terms of the generator of a two-dimensional random walk, whose jump rates are model dependent. From this, we deduce an asymptotic independence which  is shared by many  models.
\end{abstract}


%\tableofcontents

\section{Introduction}

A very challenging  problem in the field of Interacting Particle Systems (IPS) is to show its asymptotic independence. More precisely, given a sequence of continuous-time Markov chains $(\eta^N_t; t \geq 0)_{N \geq 1}$, evolving in a  discrete space $\Lambda_N$ with $N$ points and with state-space $\Omega_N$,   what is the decay, in terms of the scaling parameter $N$, of the time-dependent correlation function. Knowing this type of results is fundamental for the description of the hydrodynamic limit as well as the analysis of the fluctuations around that hydrodynamical behavior.

Over the last decades, there has been a tremendous effort to develop methods that allow approaching this problem, such as the Matrix Product Ansatz (MPA) \cite{Derrida1, Derrida2}, the {Bethe} Ansatz (BA) \cite{BAflu}, and the method of $v$-functions \cite{DeMasiTruncated,DeMasibook}. Even though the first two approaches can provide explicit formulas for stationary non-equilibrium correlations, their range of application is limited. Only special models allow for such an approach and the computation of the precise expression of these stationary correlations can be computationally complex.  On the other hand, the method of $v$-functions has been shown to be applicable for classes of models for which there is no MPA nor BA formulation, nevertheless, this approach is also applicable on very specific types of dynamics. 

We now give a brief description of the method of $v$-functions and we refer the interested reader to Section 5.4 of \cite{DeMasibook}. The study starts by analyzing the average occupation at a site $x$ and at a fixed time $t$ starting from a general probability measure $\mu^N$ defined on the state-space $\Omega_N$, i.e. \begin{equation}
\label{discrete_profile_new}
\varrho_t^N(x):=\mathbb E_{\mu^N}[\eta_t(x)].\end{equation}
A simple application of Kolmogorov's equation allows writing down the space-time evolution of this quantity as $$\partial_t\varrho_t^N(x)=\mathbb E_{\mu^N}[\mcb L_N\eta_t(x)],$$ where $\mcb L_N$ is the infinitesimal generator of our Markov process $(\eta^N_t; t \geq 0)_{N \geq 1}$. If $\mcb L\eta(x)$ is a polynomial function of $\eta's$, then the previous display might give rise to sums of terms of the form 
$\mathbb E_{\mu^N}[\eta_t(x_1)\eta_t(x_2)]$, where $x_1,x_2\in\Lambda_N$, or even higher order terms. In any case, at this point the evolution equation for $\varrho_t^N(x)$ is not closed as a function of $\varrho_t^N(x)$. What one can do then is to assume an independence of the occupation variables $\eta(x_1)$ and $\eta(x_2)$ at each time $t$, and denote the solution of the resulting equation by $\widetilde \varrho_t^N(x)$. This means that we assume that terms of the form $\mathbb E_{\mu^N}[\eta_t(x_1)\eta_t(x_2)]$ are replaced by $\mathbb E_{\mu^N}[\eta_t(x_1)]\mathbb  E_{\mu^N}[\eta_t(x_2)]$ and now we obtain an expression in terms of a function of $\widetilde \varrho_t^N$. 

Therefore, instead of looking at the $k$-points correlation function given by $$\varphi_t^N(x_1,\cdots, x_k):=
\mathbb E_{\mu^N}[(\eta_t(x_1)-\varrho_t^N(x_1))\cdots(\eta_t(x_k)-\varrho_t^N(x_k))],$$ one looks at the $v$-function defined as 
$$v_t^N(x_1,\cdots, x_k)=\mathbb E_{\mu^N}[(\eta_t(x_1)-\widetilde \varrho_t^N(x_1))\cdots(\eta_t(x_k)-\widetilde\varrho_t^N(x_k))],$$ and the goal consists in showing that 
$v_t^N(x_1,\cdots, x_k)$ decays to zero as $N\to+\infty$. 
Note that for $k=1$ the $v$-function is simply given by $v_t^N(x_1)=\varrho(x_1)-\widetilde\varrho_t^N(x_1)$ and so showing that the $v$- function vanishes with $N\to+\infty$ gives, in particular, that the solutions of the two equations are close in some sense. 

The advantage of working with $v_t^N(x_1,\cdots, x_k)$ rather than $\varphi_t^N(x_1\cdots, x_k)$ is that, for each $k \in \mathbb{N}$, the evolution equation becomes a function of known quantities, i.e. $v$-functions depending on fewer space points. In both cases, if the generator increases the degree of the polynomial functions of $\eta$, then we get in the evolution of the $k$-points correlation function a new function that depends on more than $k$ space points. As a consequence, we cannot write the equation for the $k$-points correlation function in terms of correlation functions with $k$ or fewer space points. Instead, one can iterate the replacing of higher order correlation functions and truncate this expansion at a certain point that allows controlling the error term and the remaining ones in classical ways. This program is highly technical and has proved to be suitable for many models, though the estimates can be demanding and very model-dependent.
For details we refer the reader to the book  \cite{DeMasibook}.

Here we provide an alternative approach to estimating two-points correlation functions for some classes of boundary-driven interacting particle systems, such as: the partial exclusion process, the inclusion process, {and }independent particles; interacting diffusions, such as the Ginzburg-Landau dynamics and the Brownian Energy Process;  interacting piles, such as the model introduced in \cite{GiardinaFrassek}, etc. These are just examples of models for which our machinery applies but certainly, others should be investigated. 
In due course, we will highlight what are the important features of these models that allow the application of our techniques. As the reader will see, there will be some similarity between our approach and the one to estimate $v$-functions. We observe that these two notions coincide on a variety of models as, for example, the symmetric simple exclusion process.

Inspired by previous results in the literature and in the more recent results of \cite{FGJS23}, the strategy is based on properly defining the two-points correlation function in such a way that the space-time evolution equation is written solely in terms of the correlation function itself and the function $\varrho_t^N(\cdot)$.

\subsection{Our result in a nutshell}
Let $(\eta_{tN^2})_{t \geq 0}$ be an interacting particle system with symmetric jump rates and speed up in the diffusive time scale $tN^2$ (the time scale for which there is a non-trivial evolution of the quantity of interest), such as the partial exclusion process, the inclusion process, independent particles, the Ginzburg-Landau dynamics, the Brownian Energy process and the open Harmonic model, a model of interacting piles introduced in \cite{GiardinaFrassek} and also explored in \cite{Chiara_paper}.

Let $(\mu^N)_{N \geq 1}$ be a sequence of probability measures on the state-space of the system that we denote by  $\Omega_N$ and that we write as $\Omega_N:= S^{\Lambda_N}$, with $S \subset \mathbb{R}$ and $\Lambda_N$ a discrete space with $N$ points. We remark that the set $S$ is model-dependent. We present here the results without specifying the discrete domain where particles will evolve but typically we think about the system evolving on the discrete torus $\mathbb T_N=\{0,1,\cdots, N-1\}$ where $N\equiv 0$, or on the discretization $\Lambda_N:=\{1,2,\cdots, N\}$ of the interval $[0,1]$,  but other domains could be considered. Recall \eqref{discrete_profile_new} and denote the infinitesimal generator of the process by $\mathcal L_N$.
%or on the full lattice $\mathbb Z$,
Our main result is the following theorem. 

\begin{theorem} \label{th_eq_corr}
Let $\varphi^N_t$ denote the two-points correlation function, defined on $(x,y) \in (\Lambda_N)^2$ with $y \neq x$, by
\begin{equation*}
    \varphi^N_t(x,y) = \mathbb{E}_{\mu^N} [({\eta}_{tN^2}(x)-\varrho_t^N(x)) ({\eta}_{tN^2}(y)-\varrho_t^N(y))].
\end{equation*}  
{Then,} there exists a function $f:S \to \mathbb{R}$, that can be written in terms of  the quantities $(\eta(x))^2$,$\eta(x)$ and $\varrho_t^N(x)$, and that allows defining the two-points correlation function at the diagonal points  $y=x$, i.e. for every $x \in \Lambda_N$ as
\begin{equation*}
    \varphi^N_t(x,x) = \mathbb{E}_{\mu^N}[f(\eta_{tN^2}(x))],
\end{equation*} in such a way that, under this choice, the space-time evolution of the two-points correlation functions is given  by 
\begin{equation} \label{eq_correlation_closed}
    \partial_t \varphi^N_t (x,y) = {\mcb A}_N \varphi^N_t(x,y) + g^N_t(x,y),
\end{equation} where $\mcb{A}_N$ is the infinitesimal generator of a two-dimensional  random walk with rates and state-space which are model dependent and $g^N_t$ is a function uniformly bounded in the space variables.
\end{theorem}

We observe that not only the jump rates and the state-space of the random walk, and thus the operator ${\mcb A}_N$, but also the function $g^N_t$ are all model dependent. Below we will make the derivation of this equation for a variety of models and we will see what are the corresponding random walks and the corrector function $g_t^N$. We also observe that the function $g_t^N$ appears as  a consequence of the fact that in all the interacting particle systems, the evolution of each particle is not independent of the others. In the case of the dynamics of independent particles, as expected, the function $g_t^N$ is null. 

Of course, there are many ways of rewriting the previous evolution equation in terms of different random walks. We will see below that for microscopic  models with an open boundary the scenario is exactly the same as above, i.e. the two-points correlation function also satisfies an identity as \eqref{eq_correlation_closed}. The difference only comes at the level of the random walk behavior that in this setting can be absorbed or reflected at the boundary of the domain of $\varphi_t^N$. In case the open boundary is tuned to make it fast (resp. slow) with respect to the bulk dynamics, then the random walk that we consider in the last display is the absorbed (resp. reflected) one. In the case where we write the evolution equation in terms of the reflected random walk, then it has an extra term which is of multiplicative form, but it has a negative sign and therefore does not bring any additional difficulties in deriving our results.

\subsection{Consequence of our results}
An immediate consequence of our main result is that it allows showing the decay of the two-points correlation function as $1/N$, provided that the function $f$ and the random walk with generator $\mcb A_N$ given in the theorem have suitable properties.
This is the content of the next result. 
\begin{corollary}
If \begin{equation}\label{eq:initial_corr}
    \max_{{(x,y) \in (\Lambda_N)^2}} |\varphi^N_0(x,y)|\lesssim \frac{1}{N},
\end{equation} then for every fixed $T >0$,
\begin{equation*}
    \sup_{t \in [0,T]} \max_{{(x,y) \in (\Lambda_N)^2}} |\varphi^N_t(x,y)|\lesssim \frac{1}{N}.
\end{equation*} 
\end{corollary}

The last result is saying that under the evolution equation given by \eqref{eq_correlation_closed}, the two-points correlation function as time evolves is bounded by its initial condition plus the bounds on the function $g^N_t$ (which only depends on the occupation variables $\eta(\cdot)$ i.e. polynomials of degree one in $\eta$) but also on  the occupation time of the random walk with generator $\mcb{A}_N$ (recall  \eqref{eq_correlation_closed}) on the support of $g^N_t$.

The previous corollary follows as a consequence of \eqref{eq_correlation_closed} and a simple application of Duhamel's formula plus a sharp control of the behavior of the associated random walk. Observe that  assumption \eqref{eq:initial_corr} reduces the range of possible initial distributions for our processes, nevertheless, it allows, for example,  initial product measures.

The result above represents a crucial estimate for the proof of non-equilibrium fluctuations, though it is not sufficient especially for the case of models with unbounded occupation variables, as it is the case of the inclusion process. Nevertheless, it is a step forward in showing that result for a wide range of models. We observe that for the derivation of non-equilibrium density fluctuations we typically need to control the second moments of linear functionals of the form

$$\mcb Y_t^N(f)= {\frac {1}{\sqrt{N}}}\sum_{x\in\Lambda_N}f\Big(\frac xN\Big)(\eta_{tN^2}(x)-\varrho_t^N(x)).$$

A simple computation shows that the $\mathbb L^2$ moments of this functional, with respect to $\mu^N$, can be {written as}  
\begin{equation*}
\begin{split}
\mathbb E_{\mu^N}\Big[\Big( \mcb Y_t^N(f)\Big)^2\Big]&= \frac 1 N \sum_{x\in\Lambda_N}\Big(f\Big(\frac xN \Big)\Big)^2 \mathbb E_{\mu^N}\Big[(\eta_{tN^2}(x)-\varrho_t^N(x))^2\Big]\\&+\frac{1}{N}\sum_{x\neq y\in\Lambda_N}f\Big(\frac xN \Big)f\Big(\frac yN \Big)\varphi_t^N(x,y).\\
\end{split}
\end{equation*}
So, knowing the decay as $\tfrac 1N$ of the two-points correlation function allows us to control the second term in the last identity by a term of order $O(1)$ as long as the test functions are, for example, bounded. The first term on the right-hand-side of the last display can be easily estimated if the occupation variables are bounded, but if they are not, then one needs extra arguments to control this term. 

Moreover, when the IPS is put in contact with reservoirs, knowing properties on its non-equilibrium stationary states (NESS) is a central question that needs attention.    There are techniques that have been developed over the last years to obtain information about the NESS for several models, but these are still very far from being robust techniques. 
Our approach allows to {partially} reply to this question for a variety of models. When the boundary dynamics is as strong as the bulk dynamics, one can explicitly write down the formula for the two-points stationary correlation function just by analysing the stationary solution of 
\eqref{eq_correlation_closed}, which is an harmonic function for some Laplacian operator with Dirichlet boundary conditions. When the boundary dynamics is tuned, either slow or fast, obtaining explicit expressions for the two-points stationary correlation function is hard. Nevertheless using tools from partial differential equations, such as maximum principles, we can bound the unknown solution (in the tuned case) by the known ones and obtain the desired estimates. 
As a consequence non-equilibrium stationary fluctuations, i.e. fluctuations starting from the NESS should, in principle, be possible to derived. 

\subsection{Possible future work}

In an upcoming article, we will show how to extend the results of this work to the $k$-points correlation function with $k$ being any natural number greater or equal to $2$. The generalization, which finds inspiration in \cite{Ferrari_v-functions}, is based on an application of stochastic duality and it is strongly based on the results presented in the  Appendix \ref{remark on duality} of these notes. 
As we mentioned above, knowing the decay in $N$ as $1/N$ for the two-points correlation function is necessary to control second moments of the density fluctuation field. Though this is very far from being sufficient. As in fluctuations the limit of the density fluctuation field, that we denote by $\mcb Y_t$, is a solution to some stochastic partial differential equation, one needs to argue what is the form of the noise in that equation. To do this plan, one can analyze the limit of the Dynkin's martingales associated to the density fluctuation field. This martingale has a quadratic variation which is, for some models, typically quadratic in the variables $\eta(\cdot)$. Therefore if we have estimates on the 4-points correlation function, we are then able to analyze the convergence with respect to the $\mathbb L^2$-norm of its quadratic variation.   Alternatively, what one could do it is explore the fact that the quadratic variation can be written as an additive functional of the empirical measure and then use the hydrodynamic limit to analyse the limit of the quadratic variation. This last strategy though only works if we can rewrite the quadratic variation in terms of that empirical measure and this is something that usually needs a kind of replacement lemma, which can be derived for some specific dynamics. 
An interesting application where the aforementioned routine can be done is the characterization of the non-equilibrium fluctuations of the symmetric inclusion process, which is a model with unbounded occupation variables. Nevertheless, the characterization of the noise cannot be done using the hydrodynamic limit results. This is due to the fact that the replacement lemma does not go through for this model since the partition function of its invariant measures, namely products of Negative Binomial distributions, do not satisfy the finite exponential moments' hypothesis. Such an assumption, which holds for example for the zero-range processes, is crucial in the derivation of replacement lemmas. To avoid this issue, we can then analyze the quadratic variation by controlling 4-points correlation functions. This is a subject for future work. 

We also intend to apply our routine to control  the space-time two-points correlation function given on $0\leq s\leq t\leq T$ and $x,y\in\Lambda_N$ by

\begin{equation*}
    \Psi^N_{t,x}(s,y) = \mathbb{E}_{\mu^N} [({\eta}_{tN^2}(x)-\varrho_t^N(x)) ({\eta}_{sN^2}(y)-\varrho_s^N(y))].
\end{equation*} 
The function $\Psi$ solves a much simpler equation than \eqref{eq_correlation_closed}, because the time $s$ and the space point $y$ are taken as fixed and we look at the space-time evolution of $\Psi^N_{t,x}$ as a function of $t$ and $x$. Note that its initial condition $t=s$ is given by the two-points correlation function defined above, i.e. 
$$\Psi^N_{t,x}(t,y)=\varphi^N_t(x,y).$$
Therefore,  the bounds we obtain here are also crucial to transport them to the space-time function $\Psi$. This is subject for future work. 

Another interesting program is to extend our approach in order to obtain bounds for the $k$-points (space-time) correlation functions for other interacting particle systems, especially those for which there is no duality. 
Moreover, an interesting direction to pursue would be to understand the case of asymmetric models, for which the action of the generator increases the degree of the functions. The technique shown here can not be directly applied and, to treat these types of models, further investigation will be required. More ambitiously, one could also think about extensions to the case of processes with long jumps, where, in this last case, the equation of the corresponding correlation function is governed by a fractional partial differential equation written in terms of a  non-local operator.

\subsection{Outline of the article}

In Section \ref{sub_model_def}, we introduce the models of interacting particles we will discuss; we review what is known regarding their invariant measures, and also show how to close the equation for the two-points correlation function; and we obtain bounds for such functions. In Section \ref{sub_model_def_diffusions} we approach the same type of questions and make some remarks as we did in Section \ref{sub_model_def}, but now for two types of interacting diffusions: the Ginzburg-Landau dynamics in Section \ref{GL}, and the Brownian Energy Process  in Section \ref{BEP}; and also for a model of interacting piles introduced in the literature as the (open) Harmonic Model - see Section \ref{interacting_piles_model}. Finally, in Appendix \ref{remark on duality}, we provide a remark on how to use stochastic duality to obtain directly the perturbation of the two-points correlation function for which its evolution equation can be closed.

%The correlation function is a central object in the study of both hydrodynamic limits as well as non-equilibrium fluctuations. 

\section{The models: interacting particles} 
\label{sub_model_def}

Fix a scaling parameter $N \in \mathbb N$ and let us introduce the discrete space where the particles will evolve. Let $\Lambda_N:=\{1,\dots, N-1\}$ be a discretization with {$N-1$} points of the one-dimensional interval $[0,1]$. In what follows, the set $\Lambda_N$ will always be referred to as the bulk. We say that $x,y \in \Lambda_N$ are \textit{nearest neighbors} if $|y-x|=1$, and we denote it by $x \sim y$. We introduce $\Omega_N := {S}^{\Lambda_N}$, where $S \subset \mathbb{N} \cup \{0\}$ is a given set. The set $\Omega_N$ will represent the state-space of the Markov processes that we consider, and its elements, which we denote by the greek letters $\eta$ and $\xi$, are called  \textit{configurations}. We consider Markov processes that have state-space $\Omega_N$ and that we denote by $(\eta_t)_{t \geq 0} :=(\eta_t(x); \ x\in\Lambda_N)_{t \geq 0}$, where $\eta_t(x) \in S$ denotes the quantity of particles on the configuration $\eta_t$ at the site $x\in \Lambda_N$ and at the time $t\geq 0$. 

\subsection{The dynamics} \label{sec_dynamics}

We start by defining the action of the bulk dynamics through the infinitesimal generator $\mcb{L}_{0,N}$, an operator that acts on functions $f: \Omega_N \to \mathbb R$ as
\[\label{generator_bulk}
\mcb{L}_{0,N} f(\eta) 
		:= \sum_{x=1}^{N-2} r_{x,x+1}(\eta) [ f(T_{x,x+1}(\eta)) - f(\eta)]+  r_{x+1,x}(\eta) [ f(T_{x+1,x}(\eta)) - f(\eta)],
\] 
for every $\eta \in \Omega_N$. 
Above, for each $x,y\in\Lambda_N$ the transformation $T_{x,y}:\Omega_N\to\Omega_N$ is defined as  
$$T_{x,y}(\eta):=\eta+\delta_{y}-\delta_x,$$
where $\delta_x$ is the configuration with exactly one particle which is located at the site $x
$. In this way, the transformation $T_{x,y}$ represents the removal of a particle from the site $x$ and putting it in the site $y$. This defines the bulk dynamics.

\medskip
Now we make the following assumption on the interaction rates between the sites $x$ and $y$:
\begin{equation}\label{eq:rates}
r_{x,y}(\eta)=\eta(x)(c+d\eta(y)),
\end{equation}
where $d\in \mathbb{R}$ and {$c> 0$} are such that $r_{x,y}(\eta) \geq 0$ for every $x,y \in \Lambda_N$ and every $\eta \in \Omega_N$. Now we particularize some possible choices for the rates for which we recognize some well-known processes. 

\begin{itemize}
\item
[(SEP)] Symmetric Simple Partial Exclusion Process:\\
If $c=\alpha\in\mathbb N$ and $d=-1$, then for $\Omega_N:=\{0,1,\cdots, \alpha\}^{\Lambda_N}$ the rates become equal to $$r_{x,y}(\eta)=\eta(x)(\alpha-\eta(y)),$$ and this corresponds to the partial exclusion process allowing at most $\alpha\in\mathbb N$ particles at each site.

\item [(SIP)] Symmetric Simple Inclusion Process:\\
If $c=\alpha\in\mathbb R$ and $d=1$, then for $\Omega_N:=(\mathbb{N} \cup \{0\})^{\Lambda_N}$ the rates become equal to $$r_{x,y}(\eta)=\eta(x)(\alpha+\eta(y)),$$ and this corresponds to the partial inclusion process allowing any number of particles at each site.

\item [(IRW)]  Independent Particles:\\
If  $d=0$ and $c> 0$, then for  $\Omega_N:=(\mathbb{N} \cup \{0\})^{\Lambda_N}$ the rates become equal to
 $$r_{x,y}(\eta)=c\eta(x),$$ and this corresponds to the dynamics of independent random walks. 
\end{itemize}

Now we define the boundary dynamics which is added in order to create an external interaction with the system. Under this dynamics, particles can be injected or removed in the system and this breaks down the conservation of the number of particles in the system.  The
boundary generator is an operator, denoted by $\mcb{L}_{b,N}$, that acts on functions 
 $f: \Omega_N \to \mathbb R$ as \[
\begin{split} \label{generator_boundary}
\mcb{L}_{b,N} f(\eta) 
		&:= \sum_{x\in\{0,N-1\}}r_{x,x+1}(\eta) \big\{ f(T_{x,x+1}(\eta)) - f(\eta)\big\}+  r_{x+1,x}(\eta) \big\{ f(T_{x+1,x}(\eta))- f(\eta)\big\},
\end{split}
\] 
for every $\eta \in \Omega_N$. Above 
\begin{equation}\label{eq:rates_boundary}
r_{0,1}(\eta)=\lambda_-\varrho_-(c+d\eta(1))\quad \textrm{and}\quad r_{1,0}(\eta)=\lambda_-\eta(1)(c+d\varrho_-),
\end{equation}
and, analogously,
\begin{equation}\label{eq:rates_boundary2}
r_{N-1,N}(\eta)=\lambda_+\eta(N-1)(c+d\varrho_+)\quad \textrm{and}\quad r_{N,N-1}(\eta)=\lambda_+\varrho_+(c+d\eta(N-1)),
\end{equation}
where $\lambda_\pm \in (0,1]$ and $\varrho_\pm \in (0,-\frac{c}{d})$ when $d < 0$ and $\varrho_\pm > 0$ when $d \geq 0$. Here $\varrho_-$ and $ \varrho_+$ represent the density of the left and right reservoirs, respectively. Moreover
\begin{equation}
T_{0,1}(\eta) := \eta + \delta_1 \quad \textrm{ and } \quad T_{1,0}(\eta) := \eta - \delta_1,
\end{equation} and, analogously,
\begin{equation}
T_{N,N-1}(\eta) := \eta + \delta_{N-1} \quad \textrm{ and } \quad T_{N-1,N}(\eta) := \eta - \delta_{N-1}.
\end{equation}

We consider the processes speed up in the diffusive time scale $tN^2$ in order to see a non-trivial evolution of the density of particles. In that case, the speeded process $(\eta_{tN^2})_{t \geq 0}$ has infinitesimal generator given by $$N^2\mcb L_N:=N^2\mcb L_{0,N}+N^2\mcb L_{b,N}.$$ 

The Markov process defined by the bulk generator $\mcb L_{0,N}$ conserves the total number of particles
\begin{equation*}
\sum_{x \in \Lambda_N} \eta (x).
\end{equation*}

\subsection{Invariant measures}

Here we recall what is known regarding the stationary measures for the family of processes introduced above. We divide the analysis into the following three cases.
\begin{enumerate}
\item If $d < 0$: observe that in order to have non-negative rates we must impose, for every $\eta \in \Omega_N$ and $x \in \Lambda_N$, that $c \geq -d\eta(x)$ and $\eta(x) \in \left\{0,1\,\dots,\left\lfloor -\frac{c}{d} \right\rfloor \right\}$. {Moreover, the model in this case is irreducible and with finite state-space, thus it admits a unique invariant measure.} So, we define, for every $\eta \in \Omega_N$ and $\lambda \in (0,1)$,
\begin{equation} \label{equation_inv_measure}
   \mu(\eta) := \frac{1}{Z_\mu} \prod_{x \in \Lambda_N} \frac{\lambda^{\eta(x)}(1-\lambda)^{c+d\eta(x)}}{\eta(x)! F(c+d\eta(x))},
\end{equation} where
$$F(c+d\eta(x)) :=(c+d\eta(x)) \cdot(c+d(\eta(x)+1))\cdots\left(c+d \left\lfloor -\frac{c}{d} \right\rfloor \right)$$ and $Z_\mu$ is the partition function of $\mu$,  given by $$Z_\mu := \sum_{\eta \in \Omega_N} \prod_{x \in \Lambda_N} \frac{\lambda^{\eta(x)} (1-\lambda)^{c+d\eta(x)}}{\eta(x)! F(c+d\eta(x))}.$$ If  $\lambda_- \lambda_+ (\varrho_- - \varrho_+) = 0$, taking $\lambda = \frac{\varrho_-}{c} = \frac{\varrho_+}{c}$, then  $\mu$ is the stationary measure of the model. Whenever $\lambda_- \lambda_+ (\varrho_- - \varrho_+) \neq 0$ and $c \neq 1$, we have no information on the stationary measure. For the particular case of (SEP), i.e. $c= \alpha \in \mathbb{N}$ and $d=-1$:
\begin{itemize}
\item if  $\lambda_- \lambda_+ (\varrho_- - \varrho_+) = 0$, then \eqref{equation_inv_measure} represents the product measure
\begin{align*}
\bigotimes_{x \in \Lambda_N} \textrm{Binomial}(\alpha,\lambda),
\end{align*} where $\lambda = \varrho_-/\alpha = \varrho_+/\alpha$.

\item if $\lambda_- \lambda_+ (\varrho_- - \varrho_+) \neq 0$, then for $\alpha \neq 1$ we have no information about the invariant measure. Nevertheless for $\alpha=1$ the MPA method gives information about it, for details see \cite{Derrida1}. Recently in \cite{FloreaniCasanova}, the authors have also described this invariant measure by probabilistic arguments. Though the arguments are very interesting, they do not seem easy to extend to other models as for example (SEP) under the choice $\alpha=2$. 

\end{itemize}

\item If $d > 0$: {in this case, the model is clearly irreducible and it is also possible to show that it admits in fact a unique invariant measure - see for example Appendix A of \cite{articleSFF} for the case $d=1$ and $c = \alpha > 0$, i.e. (SIP). We} define, for every $\eta \in \Omega_N$ and $\lambda \in (0,1)$,
\begin{equation} \label{inv_measure_SIP}
    \nu(\eta) := \frac{1}{Z_\nu} \prod_{x \in \Lambda_N} \frac{\lambda^{\eta(x)}(1-\lambda)^{c+d\eta(x)}}{\eta(x)! H(c+d\eta(x))},
\end{equation} where
$$H(c+d\eta(x)) := (c+d(\eta(x)-1)) \cdot (c+d(\eta(x)-2)) \cdots (c+d) \cdot c$$ and $Z_\nu$ is the partition function of $\nu$, given by $$Z_\nu := \sum_{\eta \in \Omega_N} \prod_{x \in \Lambda_N} \frac{\lambda^{\eta(x)} (1-\lambda)^{c+d\eta(x)}}{\eta(x)!H(c+d\eta(x))}.$$ If  $\lambda_- \lambda_+ (\varrho_- - \varrho_+) = 0$, taking {$\lambda = \frac{\varrho_-}{c + d\varrho_-} = \frac{\varrho_+}{c + d\varrho_+}$}, then $\nu$ is the stationary measure of the model. {If $\lambda_- \lambda_+ (\varrho_- - \varrho_+) \neq 0$, we only know that there exists a unique invariant measure for the model but we have no characterization for it.} For the particular case of (SIP), i.e. $c= \alpha \geq 0$ and $d=1$, we have that:
\begin{itemize}
    \item if  $\lambda_- \lambda_+ (\varrho_- - \varrho_+) = 0$, then \eqref{inv_measure_SIP} reduces to the product measure given by 
\begin{align*}
\bigotimes_{x \in \Lambda_N} \textrm{Negative Binomial}(\alpha,\lambda),
\end{align*}
where {$\lambda = \frac{\varrho_-}{\alpha + \varrho_-} = \frac{\varrho_+}{\alpha + \varrho_+}$}. 

    \item if $\lambda_- \lambda_+ (\varrho_- - \varrho_+) \neq 0$, we have no information on the stationary measure. 
\end{itemize}

\item If $d = 0$: it means that we are considering (IRW) and thus
\begin{itemize}
\item if $\varrho_- = \varrho_+$, then the (IRW) admits the stationary measure, which is of product form, given by
\begin{align*}
\bigotimes_{x \in \Lambda_N} \textrm{Poisson}(\lambda),
\end{align*} where $\lambda = \varrho_-/c$.

\item if $\varrho_- \neq \varrho_+$, we have no information about the invariant measures.
\end{itemize}
\end{enumerate}

\subsection{Correlation estimates} \label{subsection_correlation_IPS}

Fix a sequence of probability measures {$(\mu^N)_{N \geq 1}$} in $\Omega_N$. Let $\mathbb{P}_{\mu^N}$ be the probability measure on the Skorohod space $\mcb{D}_N([0,T],\Omega_N)$ (i.e. the space of c\`adl\`ag trajectories taking values in $\Omega_N$) induced by the Markov process $(\eta_{tN^2})_{t\geq 0}$ and by the initial measure $\mu^N$. We denote the expectation with respect to $\mathbb{P}_{\mu^N}$ by $\mathbb{E}_{\mu^N}$.

\subsubsection{The discrete density profile}

We denote the discrete density profile at time $t \geq 0$ and at site $x \in{\Lambda}_N$  by $\varrho_t^N(x)$, which is given by
\begin{equation*}
\varrho_t^N(x) := \mathbb E_{\mu^N} [ \eta_{tN^2}(x)].
\end{equation*}
First we want to understand the space-time evolution of $\varrho_t^N(x)$. To that end, 
note that taking its derivative in time we get that
$$
\partial_t \varrho_t^N(x) =\mathbb E_{\mu^N}[N^2\mcb L_N\eta_t(x)].
$$
A simple computation based on the expression of the generator gives that
$$\mcb L_{0,N}\eta(x)=r_{x-1,x}(\eta)-r_{x,x-1}(\eta)-r_{x,x+1}(\eta)+r_{x+1,x}(\eta).
$$
Since we assumed \eqref{eq:rates} we get for $x\in\{2,\cdots, N-2\}$ that 
\begin{equation}\begin{split}
\mcb L_{0,N}\eta(x)&=\eta(x-1)(c+d\eta(x))-\eta(x)(c+d\eta(x-1))\\&\quad \quad \quad \quad \quad -\eta(x)(c+d\eta(x+1))+\eta(x+1)(c+d\eta(x))\\&=c\Delta \eta(x).\end{split}\end{equation} Above we used the notation 
$\Delta \eta(x):=\eta(x-1)+\eta(x+1)-2\eta(x).$
Moreover, at the boundary, we have 
\begin{equation}\begin{split}
&\mcb L_{b,N}\eta(1)=\lambda_-\varrho_-(c+d\eta(1))-\lambda_-\eta(1)(c+d\varrho_-)=\lambda_-c(\varrho_--\eta(1)),\\
&{\mcb L_{b,N}\eta(N-1)=\lambda_+\varrho_+(c+d\eta(N-1)) -\lambda_+\eta(N-1)(c+d\varrho_+) =\lambda_+c(\varrho_+ - \eta(N-1)).}
\end{split}\end{equation}

From the last computations,   $\varrho_t^N(x)$ is the unique solution of the following initial value problem:
\begin{align} \label{system_of_eq_for_discrete_rho}
\left\{
\begin{array}{r@{\;=\;}l}
\vspace{0.1cm}
\partial_t f_t(x) & \Delta^{1D}_N f_t(x), x \in \Lambda_N, t \geq 0,\\
\vspace{0.1cm}
f_t(0) & \varrho_-, t \geq 0,\\
\vspace{0.1cm}
f_t(N) & \varrho_+, t \geq 0,\\
\vspace{0.1cm}
f_0(x) & \mathbb E_{\mu^N}[\eta_0(x)].
\end{array}
\right.
\end{align}
Above,    the operator $\Delta^{1D}_N$ is a discrete one-dimensional Laplacian defined on functions $f: \Lambda_N\cup\{0,N\} \to \mathbb R$ as
\begin{equation} \label{laplaciannn}
\Delta^{1D}_N f(x) = N^2\Big( c_{x-1,x} \big( f(x-1)-f(x) \big) + c_{x,x+1} \big( f(x+1)-f(x) \big)\Big),
\end{equation}
for every $x \in \Lambda_N$, where $c_{0,1}=c\lambda_-$ and $c_{N-1,N}=c\lambda_+$ and in all other bonds $ \{x,y\}$ the rate $c_{x,y}$ is equal to $c$.

We observe that equation \eqref{system_of_eq_for_discrete_rho} has one stationary solution. Let us denote it by $\varrho_{ss}^N(\cdot)$, which satisfies:
\begin{align} \label{eq:stat_disc_profile}
\left\{
\begin{array}{r@{\;=\;}l}
\vspace{0.1cm}
\Delta_N^{1D}\varrho_{ss}^N(x) & 0, \  x \in \Lambda_N,\\
\vspace{0.1cm}
\varrho_{ss}^N(0)  & \varrho_-,\\
\vspace{0.1cm}
\varrho_{ss}^N(N) & \varrho_+.
\end{array}
\right.
\end{align}
By direct computations,  for every $x \in \Lambda_N$ it holds
\begin{equation} \label{eq_rho_ss}
\varrho_{ss}^N(x) = a_Nx +b_N,
\end{equation} 
where
$$b_N= \frac{\lambda_+ \varrho_+(1-\lambda_-) + \lambda_-\varrho_-(1+\lambda_+(N-1))}{\lambda_+(1-\lambda_-) + \lambda_-(1+\lambda_+(N-1))}\quad \textrm{and}\quad {a_N= \frac{\lambda_+(\varrho_+ - b_N)}{1+\lambda_+(N-1)}.}$$ Observe that, if $\lambda_-=1$, then $$b_N=\varrho_-\quad \textrm{and}\quad a_N= \frac{\lambda_+(\varrho_+-\varrho_-)}{1+\lambda_+(N-1)}.$$
Also, if $\lambda_+=1$, then the expression of $a_N$ simplifies to $\frac{\varrho_+-\varrho_-}{N}$.

\subsubsection{The correlation function} \label{correlation_section}

The $k$-points correlation function indicates how the occupation variables at a given time $t$ and evaluated at 
$k$-points correlate. Here we focus on the case $k=2$ since our motivation comes from understanding the non-equilibrium density fluctuations of particle systems. If we consider our interacting particle systems with bounded variables, the knowledge of the decay in $N$ of the two-points correlation function is crucial to prove the non-equilibrium behavior of the density of particles. Nevertheless, for models with unbounded variables, we also need to know the decay of the four-points correlation function or to have knowledge on hydrodynamic limit results as we mentioned in the introduction. This comes from the fact that, typically the martingale converges to a Brownian motion, and to prove this fact, we need to extract information about the asymptotic behavior of its quadratic variation, which is a quadratic function of the given variables. The interested reader can see this calculation, for the case of (SEP), in, for example, Lemma 3.2 of \cite{FGJS23}.

To that end, we explain now an approach that allows obtaining information on the decay of the two-points correlation function in terms of the scaling parameter  $N$. As we are restricted to the case $k=2$ the correlation function is defined on the two-dimensional set  $(\Lambda_N)^2$. Moreover, since it is a symmetric function, we consider the triangle $T_N:= \{(x,y) \in (\Lambda_N)^2 \ | \ x \leq y\}$ 
as being its domain. The boundary of this set $T_N$  is given by
\[
\partial {T}_N:= \Big\{ (0,y)\,: y\in\Lambda_N\cup\{0,N\}\Big\}\cup \Big\{ (x,N)\,: x\in\Lambda_N\cup\{0,N\}\,\Big\},
\] and we denote its diagonal by
\[\mathcal{D}_N := \Big\{(x,y)\in  T_N \ | \ y = x\Big\}.\]

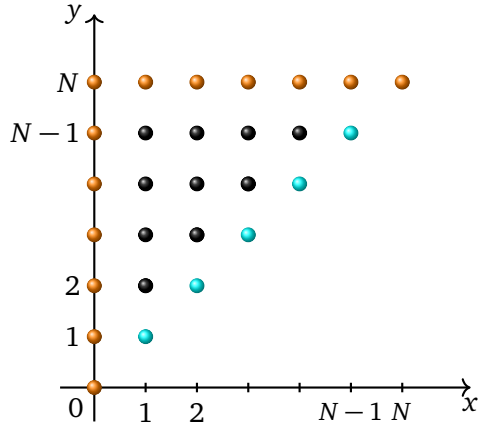
\begin{figure}[h]
    \centering
\begin{tikzpicture}[thick, scale=0.9]
\draw[->] (-0.5,0)--(5.5,0) node[anchor=north]{$x$};
\draw[->] (0,-0.5)--(0,5.5) node[anchor=east]{$y$};
\begin{scope}[scale=0.75]
\foreach \x in {1,...,3} 
	\foreach \y in {\x,...,3,4}
		\shade[ball color=black](\x,1+\y) circle (0.15);
  \foreach \x in {1,...,3,4} 
	\foreach \y in {\x}
		\shade[ball color=black](\x,1+\y) circle (0.15);
\foreach \x in {1,...,5} 
	\shade[ball color = orange](\x,6) circle (0.15); 
\foreach \x in {1,...,5} 
	\shade[ball color = cyan](\x,\x) circle (0.15); 
\foreach \y in {1,...,5} 
	\shade[ball color = orange](0,\y) circle (0.15); 
%%%%%%%
%%%%%%%
\shade[ball color= orange](6,6) circle (0.15);
\shade[ball color= orange](0,6) circle (0.15);
\shade[ball color= orange](0,0) circle (0.15);
\end{scope}	
\draw (0,0) node[anchor=north east] {$0$};
\draw (0.75,2pt)--(0.75,-2pt) node[anchor=north] {$1$};
\draw (1.5,2pt)--(1.5,-2pt) node[anchor=north]{$2$};
\draw (2.25,2pt)--(2.25,-2pt);
\draw (3,2pt)--(3,-2pt);
\draw (3.75,2pt)--(3.75,-2pt) node[anchor=north]{\small $N-1$};
\draw (4.5,2pt)--(4.5,-2pt) node[anchor=north]{\small $N$};
\draw (-0.05,.75) node[anchor=east] {$1$};
\draw (-0.05,1.5) node[anchor=east]{$2$};
\draw (-0.05,3.75) node[anchor= east]{$N-1$};
\draw (-0.05,4.5) node[anchor=east]{$N$};
\end{tikzpicture}
\caption{Illustration of the sets  $\partial T_N$ (in orange) and $\mathcal{D}_N$ (in cyan).}
    \label{fig:triangle1}
\end{figure}
%%%%%% End of triangle %%%%%%%%%%%%

\begin{definition}[Two-points correlation function]
\quad

Let  $\varphi_t^N$ be the two-points correlation function at time $t$, which is defined on $(x,y) \in T_N$ with $x \neq y$ by
\begin{equation} \label{time_dependent_correlation_def}
 \varphi^N_t(x,y)  =
\mathbb{E}_{\mu^N} [\overline{\eta}_{tN^2} (x)\overline{\eta}_{tN^2} (y)], \textrm{ if }  (x,y) \notin \partial T_N
\end{equation}
and at the boundary $\partial T_N$ we prescribe it to be null. Above $\overline{\eta} (x) := \eta(x) - \mathbb{E}_{\mu^N} [\eta(x)].$
\end{definition}

\begin{remark}
Observe that when the occupation variables take only the values $0$ or $1$ as in the (SEP) with $\alpha=1$, then the correlation function does not need to be defined at the diagonal line $\mathcal{D}_N$, since in (SEP) with $\alpha=1$ particles do not overlap.
\end{remark}

Up to this point, we have not defined   what should be the value of the correlation function $\varphi^N_t$ when extended to points in $\mathcal{D}_N$. In order to do that we will first try to find an evolution equation for $\varphi^N_t(x,y)$ that involves the generator of a two-dimensional random walk and choose the value of $\varphi^N_t$ at the points of $\mathcal{D}_N$ in such a way that we force the evolution equation to be closed. With this in mind, we note that
\begin{equation}\label{eq:evol_eq_corr}
\partial_t \varphi^N_t(x,y)=\partial_t\mathbb{E}_{\mu^N} [{\eta}_{tN^2}(x){\eta}_{tN^2}(y)]-\partial_t(\varrho^N_t(x)\varrho^N_t(y)).
\end{equation}
From \eqref{system_of_eq_for_discrete_rho} we can compute the rightmost term on the right-hand side of the last display and so it remains to analyze $\partial_t\mathbb{E}_{\mu^N} [{\eta}_{tN^2}(x){\eta}_{tN^2}(y)]$.
Using Chapman-Kolmogorov's equation, it holds 
\begin{equation}
\partial_t\mathbb{E}_{\mu^N} \Big[{\eta}_{tN^2}(x){\eta}_{tN^2}(y)\Big]=\mathbb{E}_{\mu^N} \Big [N^2\mcb L_N({\eta}_{tN^2} (x){\eta}_{tN^2}(y))\Big  ].
\end{equation}

We have now to compute the action of $\mcb L_N$ over the function $f(\eta) = \eta(x)\eta(y)$. We concentrate on the contribution of the bulk dynamics since it is there that all the difficulty arises and we will explain in due course what are the main problems in closing the evolution equation for the correlation function when boundary dynamics are present. Note that, for $|x-y|>1$, we have that
\begin{equation*}\begin{split}
\mcb L_{0,N}(\eta(x)\eta(y))=& (r_{x-1,x}-r_{x,x-1})\eta(y)- (r_{x,x+1}-r_{x+1,x})\eta(y)\\+& (r_{y-1,y}-r_{y,y-1})\eta(x)- (r_{y,y+1}-r_{y+1,y})\eta(x)~\\=&
\eta(y)\mcb L_{0,N}(\eta(x)) +\eta(x)\mcb L_{0,N}(\eta(y)),
\end{split}
\end{equation*}
and thus
\begin{equation*}
\begin{split}
N^2 \mcb L_{0,N}(\eta(x)\eta(y))= N^2\eta(y)c\Delta \eta(x) +N^2\eta(x)c\Delta \eta (y),
\end{split}
\end{equation*}
which, together with \eqref{system_of_eq_for_discrete_rho}, implies that
\begin{equation*}
\partial_t \varphi^N_t(x,y)=\Delta^{2D}_N\varphi^N_t(x,y),
\end{equation*}
where $\Delta^{2D}_N$ is the  two-dimensional discrete Laplacian which is the operator defined, for every $f:T_N\to\mathbb R$ and $(z,w) \in T_N$ with $|z-w| \geq 2$, $z \neq 1$ and $w \neq N-1$, by
\begin{align*}
    \Delta^{2D}_Nf(z,w) = cN^2(f(z+1,w) + f(z-1,w) + f(z,w+1) + f(z,w-1) - 4 f(z,w)).
\end{align*}
Now we look at neighboring points i.e. points laying in the upper diagonal  $$\mathcal{D}_N^+:= \Big\{(x,y)\in T_N \ | \ y=x + 1\Big\}.$$

Observe that 
\begin{align*}\nonumber
\mcb L_{0,N}(\eta(x)\eta(x+1))=&( r_{x-1,x}-r_{x,x-1})\eta(x+1)+ (-r_{x+1,x+2}+r_{x+2,x+1})\eta(x)\\ \nonumber
+& r_{x,x+1}\Big\{(\eta(x)-1)(\eta(x+1)+1)-\eta(x)\eta(x+1) \Big\}
\\ \nonumber +& {r_{x+1,x}\Big\{(\eta(x)+1)(\eta(x+1)-1)-\eta(x)\eta(x+1) \Big\}}
\\ \nonumber =&c( \eta(x-1)-\eta(x))\eta(x+1)- c(\eta(x+1)-\eta(x+2))\eta(x)\\ 
+&c\eta(x)^2-c\eta(x)-2(c+d)\eta(x)\eta(x+1)+c\eta(x+1)^2-c\eta(x+1).
\end{align*}
{Remark that, if $\eta(x) \in \{0,1\}$ and $c+d=0$, i.e. in the case of (SEP) with $\alpha=1$, then the sum of all the terms in the  last line of last display is equal to zero.}
If we now combine {the previous computation} with the contribution from the rightmost term in \eqref{eq:evol_eq_corr} we get
\begin{equation*}\begin{split}
\partial_t \varphi^N_t(x,x+1)&=c N^2\mathbb E_{\mu^N}[\eta_{tN^2}(x-1)\eta_{tN^2}(x+1)]-2 N^2\mathbb E_{\mu^N}[\eta_{tN^2}(x)\eta_{tN^2}(x+1)]\\
&+c N^2\mathbb E_{\mu^N}[\eta_{tN^2}(x)\eta_{tN^2} (x+2)] +c N^2\mathbb E_{\mu^N}[\eta_{tN^2}^2(x)]-cN^2\varrho_t^N(x)\\
&-2(c+d)N^2\mathbb E_{\mu^N}[\eta_{tN^2}(x)\eta_{tN^2} (x+1)]+cN^2\mathbb E_{\mu^N}[\eta_{tN^2}^2(x+1)]-cN^2\varrho_t^N(x+1)\\&-\varrho_t^N(x+1)\Delta^{1D}_N \varrho_t^N(x) -\varrho_t^N(x)\Delta^{1D}_N \varrho_t^N(x+1). 
\end{split}
\end{equation*}
By rewriting the terms in the last display in terms of the correlation function we get
\begin{equation}\begin{split}\label{eq:rates_nei}
\partial_t \varphi^N_t(x,x+1)=&cN^2\varphi^N_t(x-1,x+1)-2N^2\varphi^N_t(x,x+1)+cN^2\varphi^N_t(x,x+2)\\
+&cN^2\mathbb E_{\mu^N}[\eta_{tN^2}^2(x)]-cN^2\varrho_t^N(x)-cN^2\varrho_t^N(x)^2\\+&cN^2\mathbb E_{\mu^N}[\eta_{tN^2}^2(x+1)]-cN^2\varrho_t^N(x+1)-cN^2\rho_t^N(x+1)^2\\-&2d N^2 \varrho_t^N(x)\varrho_t^N(x+1) -2(c+d) N^2\varphi_t^N(x,x+1).
\end{split}
\end{equation}
We observe that the last identity is a consequence of our particular choice of the jump rates since we did not get terms of higher order correlation functions or even more complicated functions of the occupation variables $\eta(x)$ as it typically appears in,  for example, asymmetric models. Now, if we denote $$(c+d)G(x) :=c\mathbb E_{\mu^N}[\eta_{tN^2}^2(x)]-c\varrho_t^N(x)-(c+d)\varrho_t^N(x)^2,$$
we see that \eqref{eq:rates_nei} rewrites as
\begin{equation*}\begin{split}
\partial_t \varphi^N_t(x,x+1)=&c N^2\varphi^N_t(x-1,x+1)-2 N^2\varphi^N_t(x,x+1)+c N^2\varphi^N_t(x,x+2)\\
+&(c+d) N^2[G(x)+ G(x+1)- 2\varphi_t^N(x,x+1)]+d N^2(\varrho_t^N(x)-\varrho_t^N(x+1))^2.
\end{split}
\end{equation*}

This somehow forces {us, whenever $c+d \neq 0$,} to define the correlation function at the diagonal $\mathcal D_N$ as being equal to $G(\cdot)$ i.e. $$\varphi_t^N(x,x):=\mathbb E_{\mu^N}\Big[\frac{c}{c+d}\eta_{tN^2}(x)(\eta_{tN^2}(x)-1)\Big]-\varrho_t^N(x)^2
.$$ 
We remark that in order to obtain the previous definition, we could have also used some duality arguments - see Appendix \ref{remark on duality} for connections of the definition of the function $\varphi^N_t$  with the notion of stochastic duality.
Finally, given this definition at the diagonal {whenever $c+d \neq 0$} and noting that
\begin{equation*}\begin{split}
\mcb L_{0,N}{\eta}^2(x)=&(r_{x-1,x}+r_{x+1,x})((\eta(x)+1)^2)-\eta(x)^2)+(r_{x,x-1}+r_{x,x+1})((\eta(x)-1)^2-\eta(x)^2)\\=&2(c+d)\eta(x)[\eta(x-1) + \eta(x+1)]+c[\eta(x-1)+\eta(x+1)]-4c \eta^2(x)+2c \eta(x),
\end{split}
\end{equation*}
we see that
\begin{equation*}\begin{split}
\partial_t \varphi^N_t(x,x) =&\partial_t \mathbb E_{\mu^N}\Big[\frac{c}{c+d}\eta_{tN^2}(x)(\eta_{tN^2}(x)-1)\Big]-2\varrho_t^N(x)\Delta^{1D}_N\varrho_t^N(x)
\\
=&2cN^2\varphi^N_t(x-1,x)+2cN^2 \varphi_t^N(x,x+1)-4cN^2\varphi_t^N(x,x).
\end{split}
\end{equation*}
At this point, we have that
\begin{equation} \label{eq_varphi}
\partial_t \varphi_t^N(x,y)= \Delta^{2D}_N\varphi_t^N(x,y)+d N^2(\varrho_t^N(x+1) - \varrho_t^N(x))^2\delta_{\{y=x+1\}}(x,y),
\end{equation}
with $\Delta^{2D}_N$ being the generator of a two dimensional  random walk denoted by $(\mathcal  X_{tN^2})_{ t \geq 0}$ that  evolves on the triangle $T_N$ and moves to nearest-neighbors  at rate $c$, except at $\mathcal D_N^+$ where it moves right/up at rate $c$ and left/down at rate $(c+d)$ and it is reflected at  $\mathcal D_N^+$ if the occupation variables take values in $\{0,1\}$, otherwise it is reflected at $\mathcal{D}_N$. Moreover, repeating the computation above but with respect to the boundary dynamics, we see that it is absorbed at $\partial T_N$ with rate $c\lambda_-$ at the set of points $$\Big\{(0,y): y=0,\cdots, N\Big\}$$ and with rate $c\lambda_+$ at the set of points $$\Big\{(x,N): x=0,\cdots, N\Big\}.$$ Thus, $\Delta^{2D}_N$ is the operator acting on $f:T_N \cup \partial T_N\to\mathbb R$ such that $f(x,y) = 0$ for every $(x,y) \in \partial T_N$, as
\begin{equation} \begin{split} \label{dual_def} 
(\Delta^{2D}_N f)(x,y) &=  N^2(c_{x,x-1}[f(x-1,y)-f(x,y)] + c_{x,x+1}[f(x+1,y)-f(x,y)])\\
&+  N^2(c_{y,y-1}[f(x,y-1)-f(x,y)] + c_{y,y+1}[f(x,y+1)-f(x,y)])\\ 
&+ dN^2[f(x,x) + f(x+1,x+1)- 2f(x,x+1)]\mathbb{1}((x,y) \in \mathcal{D}^+_N),
\end{split}
\end{equation} where, for every $z \in \Lambda_N$,
\begin{align} \label{ratesss}
c_{z,z-1} = \begin{cases} 
c \lambda_-, \textrm{ if } z = 1,\\
c, \textrm{ if } z \geq 2,
\end{cases} \quad \textrm{ and } \quad 
c_{z,z+1} = \begin{cases} 
c \lambda_+, \textrm{ if } z = N-1,\\
c, \textrm{ if } z \leq N-2.
\end{cases}
\end{align}

\begin{figure}[ht]
\centering
\begin{tikzpicture}[thick, scale=1]
\draw[->] (-0.5,0)--(5.5,0) node[anchor=north]{$x$};
\draw[->] (0,-0.5)--(0,5.5) node[anchor=east]{$y$};
\begin{scope}[scale=0.75]
\foreach \x in {1,...,3} 
	\foreach \y in {\x,...,3,4}
		\shade[ball color=black](\x,1+\y) circle (0.15);
  \foreach \x in {1,...,3,4} 
	\foreach \y in {\x}
		\shade[ball color=black](\x,1+\y) circle (0.15);
\foreach \x in {1,...,5} 
	\shade[ball color =  orange](\x,6) circle (0.15); 
\foreach \x in {1,...,5} 
	\shade[ball color = cyan](\x,\x) circle (0.15); 
\foreach \y in {1,...,5} 
	\shade[ball color =  orange](0,\y) circle (0.15); 
%%%%%%%
\shade[ball color=  orange](6,6) circle (0.15);
\shade[ball color=  orange](0,6) circle (0.15);
\shade[ball color=  orange](0,0) circle (0.15);
\end{scope}	
\draw (0,0) node[anchor=north east] {$0$};
\draw (0.75,2pt)--(0.75,-2pt) node[anchor=north] {$1$};
\draw (1.5,2pt)--(1.5,-2pt) node[anchor=north]{$2$};
\draw (2.25,2pt)--(2.25,-2pt);
\draw (3,2pt)--(3,-2pt);
\draw (3.75,2pt)--(3.75,-2pt) node[anchor=north]{\small $N-1$};
\draw (4.5,2pt)--(4.5,-2pt) node[anchor=north]{\small $N$};
\draw (-0.05,.75) node[anchor=east] {$1$};
\draw (-0.05,1.5) node[anchor=east]{$2$};
\draw (-0.05,3.75) node[anchor= east]{$N-1$};
\draw (-0.05,4.5) node[anchor=east]{$N$};
%%%%%% Arrows %%%%%%%%
\node at (1.5,3) (C) { };
\node at (0.75,3) (Cl) { };
\node at (2.25,3) (Cr) { };
\node at (1.5,3.75) (Cup) { };
\node at (1.5,2.25) (Cdow) { };
\node at (0.75,3.75) (B) { };
\node at (0,3.75) (Bl) { };
\node at (0.75,4.5) (Bup) { };
\node at (3,3.75) (D) { };
\node at (2.25,3.75) (Dl) { };
\node at (3,4.5) (Dup) { };
%%%%% Arrows diag %%%%%%%%
\node at (3,3.75) (Di2) { };
\node at (3,3) (Dow2) { };
\node at (1.5,2.25) (Di) { };
\node at (0.75,2.25) (Dil) { };
\node at (2.25,2.25) (Dir) { };
\node at (1.5,3) (Diup) { };
\node at (1.5,1.5) (Dow) { };
\draw[-latex] (Dow2) to[out=60,in=-60] node[midway,font=\scriptsize,right] {$2c$} (Di2);
\draw[-latex] (Dow2) to[out=-120,in=-60] node[midway,font=\scriptsize,below] {$2c$} (Cr);
\draw[-latex] (Di) to[out=-60,in=-120] node[midway,font=\scriptsize,below,color=blue] {$c+d$} (Dir);
\draw[-latex] (Di) to[out=-130,in=130] node[midway,font=\scriptsize,left,,color=blue] {$c+d$} (Dow);
\draw[-latex] (Di2) to[out=-130,in=130] node[midway,font=\scriptsize,left,color=blue] {$c+d$} (Dow2);
%%%% end arrows diag %%%%%%%
\draw[-latex] (C) to[out=130,in=60] node[midway,font=\scriptsize,above] {$c$} (Cl);
\draw[-latex] (C) to[out=-60,in=-120] node[midway,font=\scriptsize,below] {$c$} (Cr);
\draw[-latex] (C) to[out=60,in=-60] node[midway,font=\scriptsize,right] {$c$} (Cup);
\draw[-latex] (C) to[out=-130,in=130] node[midway,font=\scriptsize,left] {$c$} (Cdow);
\draw[-latex] (B) to[out=130,in=60] node[midway,font=\scriptsize,above] {\textcolor{orange}{$c \lambda_-$}} (Bl);
\draw[-latex] (D) to[out=130,in=60] node[midway,font=\scriptsize,above] {$c$} (Dl);
\draw[-latex] (D) to[out=60,in=-60] node[midway,font=\scriptsize,right] {\textcolor{orange}{$c \lambda_+$}} (Dup);
\end{tikzpicture}
\caption{Jump rates of the random walk $(\mathcal  X_{t})_{t \geq 0}$.}
\label{fig:triangle}
\end{figure}
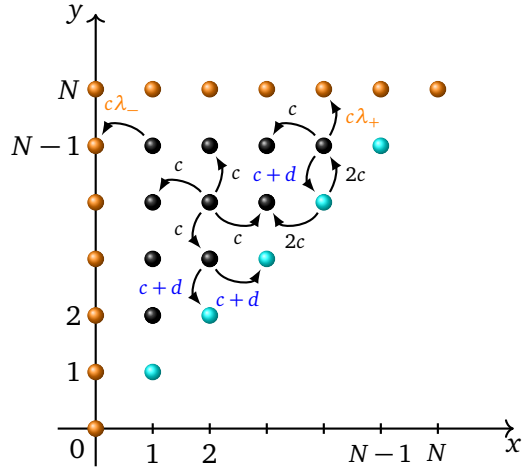

Figure \ref{fig:triangle} is an illustration of the jump rates of $(\mathcal  X_{t})_{t \geq 0}$ which are associated with the different coefficients of the operator $\Delta^{2D}_N$ without the time speed $N^2$.

%%%%%%%% Picture with jump rates %%%%%%%%%

\begin{remark}
Observe that in the case of (SEP) with $c=\alpha=1$, which means that at most  one particle is allowed at each site, the random walk jumps to the diagonal with a rate equal to $c+d = 0$.  Since in this case $c=1,d=-1$ this means that jumps to the diagonal  are suppressed. Moreover, for the (IRW), since $d=0$, the random walk jumps to the diagonal $\mathcal{D}_N$ with a rate equal to $c$. In this case, the identity \eqref{eq_varphi} simplifies to 
$$\partial_t \varphi_t^N(x,y)=\Delta^{2D}_N\varphi_t^N(x,y).$$
\end{remark}
Since \eqref{eq_varphi} holds, then, from Duhamel's formula, we know that
\begin{align} \label{varphi_SEP_1_rep}
    \varphi_t^N(x,y) = \mathbb{E}_{(x,y)} \Big[\varphi_0^N(\mathcal  X_{tN^2})+ \int_0^t g^{N}_{t-s}( \mathcal  X_{sN^2})\mathbb {1}( \mathcal  X_{sN^2}\in \mathcal{D}_N^+) \,ds\Big],
\end{align}where $\mathbb{E}_{(x,y)}$ denotes the expectation with respect to the law of the random walk $(\mathcal  X_{t})_{t \geq 0}$ starting from $(x,y)$. Above 
$g^{N}_t(z,w) =d N^2(\varrho_t^N(x+1) - \varrho_t^N(x))^2\delta_{\{y=x+1\}}(z,w).
$

Let us now focus on the rightmost term in \eqref{varphi_SEP_1_rep}. Note that 
{\begin{align*} 
\int_0^t g^{N}_{t-s}( \mathcal  X_{sN^2})\mathbb {1}( \mathcal  X_{sN^2} \in \mathcal{D}_N^+) \,ds \leq  \sup_{r\in{[0,T]}}\max_{x\in \Lambda_{N-1}}|g_r^N(x,x+1)|\int_0^t \mathbb {1}( \mathcal  X_{sN^2}\in\mathcal D_N^+)ds.
\end{align*}}
We can then take the time integral up to $t=\infty$ to get that
\begin{align}
\label{chacabuco}
{\max_{\substack{(x,y) \in T_N}}  |\varphi_t^N(x,y)| 
    		\leq  \max_{\substack{(z,w) \in T_N}}  |\varphi_0^N(z,w)| + \sup_{r \geq 0}  \max_{z \in \Lambda_{N-1}} |g_r^N(z,z+1)| \max_{\substack{(x,y) \in T_N}}     \mathcal{T}_N(x,y),}
\end{align} where 
\begin{equation}
\label{chena}
    \mathcal{T}_N(x,y) := \mathbb{E}_{(x,y)} \left[\int_0^\infty \mathbb{1}(\mathcal X_{tN^2} \in \mathcal{D}_N^+) dt \right]
\end{equation}
is the expected occupation time of the diagonal $\mathcal{D}_N^+$ by the random walk $(\mathcal X_{tN^2})_{t \geq 0}$. 

We observe that, if $d=0$, then $g^N_r = 0$ for all $r \in [0,T]$. Thus, in order to estimate $${\max_{\substack{(x,y) \in T_N}}}  |\varphi_t^N(x,y)|,$$ it would be enough to estimate ${\max_{\substack{(x,y) \in T_N}}  |\varphi_0^N(x,y)|}$. Since this only depends on the initial sequence of measures $(\mu^N)_{N \in \mathbb{N}}$, this case is done. We are only left to treat the case $d \neq 0$. 

Since for $d \neq 0$ we have $g^N_t \neq 0$, recalling \eqref{chacabuco}, we need now to provide a bound in $N$ for ${\max_{\substack{(x,y) \in T_N}}} \mathcal{T}_N(x,y)$. A simple consequence of the Markov property is that $\mathcal{T}_N$ is solution to
$$\Delta^{2D}_N \mathcal{T}_N(x,y)=-\delta_{\{y=x+1\}}(x,y),$$
and so can be explicitly computed, at least for the choice $\lambda_-=\lambda_+=1$. In this case, 
\begin{equation}\label{ot_equal_parameters}
\mathcal{T}_N(x,y) = \frac{1}{N^2}\frac{(N-y)x}{cN+d}-\frac{1}{2N(cN+d)}\mathbb{1}_{x=y},
\end{equation} which shows that, if $\lambda_-=\lambda_+=1$, then
\begin{equation} \label{cota_occupation_time_in_N}
    {\max_{\substack{(x,y) \in T_N}}}     \mathcal{T}_N(x,y) = O\left(\frac{1}{N}\right).
\end{equation}
If we do not impose $\lambda_-=\lambda_+=1$, the solution is not of polynomial form. Nevertheless, a simple application of the discrete maximum principle - see  Theorem A.1 of \cite{FGJS23} - also shows that \eqref{cota_occupation_time_in_N} still holds when taking $\lambda_-$ and $\lambda_-$ strictly positive but not necessarily equal (and in particular not both equal to $1$). In fact, we note that
\begin{equation*}
\begin{split}
\Delta^{2D}_N \Big(\mathcal{T}^{\lambda_-,\lambda_+}_N(x,y) - \mathcal{T}^{1,1}_N(x,y)\Big)&=-\alpha N^2 \left[(\lambda_- - 1)[1-\mathbb{1}_{y=1}]\mathcal{T}^{1,1}_N(1,y) \mathbb{1}_{x = 1}\right. \\&\left.+(\lambda_+ - 1)[1-\mathbb{1}_{x=N-1}]\mathcal{T}^{1,1}_N(x,N-1) \mathbb{1}_{y = N-1}\right],
\end{split}
\end{equation*}
where $\mathcal{T}^{\lambda_-,\lambda_+}_N$ (resp. $\mathcal{T}^{1,1}_N$) represents the occupation time of the random walk  on $\mathcal D_N^+$ when taking in the boundary rates  $\lambda_-,\lambda_+ \in (0,1]$ 
(resp.  $\lambda_- = \lambda_+ =1$). 
Note that, since $\mathcal{T}^{1,1}_N \geq 0$ and $\lambda_-,\lambda_+ \in (0,1]$, we conclude that the last display is greater or equal to zero. Thus, from Theorem A.1 of \cite{FGJS23}, $\mathcal{T}_N(x,y) \leq \mathcal{T}^{1,1}_N(x,y)$, for every $(x,y) \in T_N$.

\paragraph{Stationary correlations}
We observe that the stationary solution of\eqref{eq_varphi}, i.e. the solution of 
\begin{equation*}
\Delta^{2D}_N\varphi_{ss}^N(x,y)=-dN^2(\varrho_{ss}^N(x+1) - \varrho_{ss}^N(x))^2\delta_{\{y=x+1\}}(x,y),
\end{equation*}
is given, when $c=1$ and $d=-1$ for every $(x,y) \in T_N$ with $y \neq x$, by
\begin{equation} \label{relation_occupation_time_stationary_corr}
\begin{split}
\varphi_{ss}^N(x,y)=d N^2(\varrho_{ss}^N(x+1) - \varrho_{ss}^N(x))^2\mathcal{T}_N(x,y),
\end{split}
\end{equation} {and the same expression holds when $c\neq1$ or $d\neq-1$ for every $(x,y) \in T_N$. Taking }$\lambda_-=\lambda_+=1$, the function $\varphi_{ss}^N$ is explicitly given by
\begin{equation} \label{corr_ss}
\begin{split}
\varphi_{ss}^N(x,y)=d(\varrho_+-\varrho_{-})^2\Big(\frac{1}{N^2}\frac{(N-y)x}{cN+d}-\frac{1}{2N(cN+d)}\mathbb {1}_{x=y}\Big),
\end{split}
\end{equation} where here we used the explicit expression for $\varrho_{ss}^N(x)$ given in \eqref{eq_rho_ss}.
We recall that when $c=1$ and $d=-1$, we do not define $\varphi^N_{ss}(x,x)$. Finally, since  \eqref{relation_occupation_time_stationary_corr} holds, using the expression of $\varrho^N_{ss}$ given in \eqref{eq_rho_ss}, we also have that
\begin{equation*}
    \max_{\substack{(x,y) \in T_N }}     \varphi^N_{ss}(x,y) = O\left(\frac{1}{N}\right).
\end{equation*}

\section{Other models} 
\label{sub_model_def_diffusions}
In this section we provide some examples of microscopic models whose variables can take real values. In the next subsection, we present the strategy of our estimates for some examples as  the Ginzburg-Landau dynamics, the Brownian energy process and the Harmonic model. 
\subsection{Interacting diffusions}
\subsubsection{Ginzburg-Landau dynamics} \label{GL}

Recall the definition of $\Lambda_N=\{1, \ldots,N-1\}$. Fix two real constants $\Phi_-$ and $\Phi_+$.
The boundary-driven Ginzburg-Landau dynamics with linear interaction evolving on the discrete lattice $\Lambda_N$ is the Markov process $(\varphi_t)_{t \geq 0} := (\varphi_t (x) \, : \, x \in \Lambda _N)_{ t\geq 0 }$ with state-space $\Omega_N=\mathbb R^{ \Lambda_N}$ and whose generator $\widetilde {\mcb L}_N$ is given by
\begin{equation*}
N^2\widetilde {\mcb L}_N =N^2\widetilde {\mcb L}_- +N^2\widetilde{\mcb L}_+ +N^2\widetilde{\mcb L}_{b,N}.
\end{equation*}
The bulk generator is given by
\begin{equation}
\label{eq:Lb}
\begin{split}
\widetilde{\mcb L}_{b,N} &= \sum_{x \in \Lambda_{N-1}}   \left\{  ({\varphi (x)} -{\varphi (x+ 1)} ) \  (\partial_{\varphi (x+ 1)} -\partial_{\varphi (x)} )  +  (\partial_{\varphi (x+ 1)} -\partial_{\varphi (x)} )^2 \right\},
\end{split}
\end{equation}
and the boundary generators by 
\begin{equation}
\label{eq:Lbo}
\widetilde{\mcb L}_- = \partial^2_{\varphi(1)} + (\Phi_- -\varphi (1))\partial_{\varphi(1)}, \quad  \widetilde{\mcb L}_+ = \partial^2_{\varphi(N-1)} + (\Phi_+ -\varphi (N-1))\partial_{\varphi(N-1)}.
\end{equation}
In this model the quantity 
\begin{equation*}
\mathcal E(\varphi):=\cfrac12 \sum_{x \in \Lambda_N} \varphi^2 (x) 
\end{equation*}
denotes the energy of the configuration $\varphi$. 
Moreover, the Markov process defined by the bulk generator $\widetilde{\mcb L}_{N,b}$ conserves the volume 
\begin{equation*}
\mathcal V(\varphi):=\sum_{x \in \Lambda_N} \varphi (x) 
\end{equation*}
of the interface. 

\paragraph{Invariant measures}

When  $\Phi_-=\Phi_+=\Phi$, the Gaussian (product) probability measure $\mu_\Phi$ given by 
\begin{equation*}
\mu_\Phi (d\varphi) = {\mathcal Z}_\Phi^{-1} \exp\left( -\mathcal E (\varphi) +\Phi \mathcal V (\varphi) \right) d\varphi
\end{equation*}
is the unique invariant probability measure of the dynamics. Here $\mathcal Z_\Phi = Z^{N-1} (\Phi)$ where
 \begin{equation*}
Z(\Phi) =\sqrt{2\pi} \ e^{\Phi^2 / 2}.
\end{equation*} Moreover, $\mu_\Phi$ is in fact reversible. 
To prove this last fact is is enough to compute 
$\int \widetilde{\mcb L} f(\varphi) g(\varphi)d\mu_{\Phi}$ for any $f,g$ and to use the fact that the measure $\mu_{\Phi}$ is product and that the generator $\widetilde {\mcb L}$ only changes the variables in at most two space points. Then we make two summation by parts and use the explicit expression for the invariant measure to write the last identity as $\int f(\varphi) \widetilde {\mcb L} g(\varphi)d\mu_{\Phi}$. We leave the details to the reader.  As in the case of particle models, when the density of the reservoirs is such that $\Phi_-\neq \Phi_+$, we do not have any information on the invariant measure $\mu_{ss}$. 

\paragraph{The correlation function}

Let {$(\mu^N)_{N \geq 1}$} be a sequence of probability measures in $\Omega_N$. Following the same type of strategy of Section \ref{correlation_section}, we define for every $x \in \Lambda_N$,
$$v^N_t(x) := \mathbb{E}_{\mu^N}[\varphi_{tN^2}(x)],$$ and also define the two-points correlation function $\psi^N_t$ for the Ginzburg-Landau dynamics as
\begin{equation}
\psi^N_t(x,y) := \mathbb{E}_{\mu^N}[\varphi_{tN^2}(x) \varphi_{tN^2}(y)] - \mathbb{E}_{\mu^N}[\varphi_{tN^2}(x)] \mathbb{E}_{\mu^N}[\varphi_{tN^2}(y)],
\end{equation} for every $(x,y) \in T_N$ and, for $(x,y) \in \partial T_N$ we set $\psi^N_t (x,y) = 0$. A simple computation shows that
\begin{align*}
\partial_t v^N_t(x)  = \Delta^{1D}_N v^N_t(x),
\end{align*} where $\Delta^{1D}_N$ coincides with the operator \eqref{laplaciannn} by taking the rates $c_{x,x+1}=1$ for all $x$ and defining $v^N_t(0)=\Phi_-$ and $v^N_t(N)=\Phi_+$. Moreover, for every $(x,y) \in T_N$,
\begin{align} \label{equation_corr_GL}
\partial_t \psi^N_t(x,y) = \Delta^{2D}_N \psi^N_t(x,y) - 2N^2 \mathbb{1}_{y=x+1} + 4 N^2 \mathbb{1}_{y=x},
\end{align} where $\Delta^{2D}_N$ corresponds to the operator given in \eqref{dual_def} with $d=0$ and $c_{x,x+1}=c_{x+1,x}=1$.

Observe now that even though we can write the solution of \eqref{equation_corr_GL} with a stochastic representation, it is no longer true that
\begin{equation*}
    - 2N^2 \mathbb{1}_{y=x+1} + 4 N^2 \mathbb{1}_{y=x}
\end{equation*} is bounded by a constant uniformly in $N$. But, due to the form of the equation \eqref{equation_corr_GL}, we can still fix this problem in order to obtain a good bound in $N$ for $\sup_{t \in[0,T]} {\max_{\substack{(x,y) \in T_N}} |\psi^N_t(x,y)|}$. Set
\begin{equation}
    \widetilde{\psi}^N_t(x,y) := \psi^N_t(x,y) - \mathbb{1}_{y=x}.
\end{equation} Then, a simple computation shows that
\begin{equation} \label{new_eq_corr_GL}
    \partial_t \widetilde{\psi}^N_t(x,y) =  \Delta^{2D}_N \widetilde{\psi}^N_t(x,y).
\end{equation} Therefore
\begin{equation*}
    \widetilde{\psi}^N_t(x,y) = \mathbb{E}_{(x,y)}\left[\widetilde{\psi}^N_0 \left(\mathcal{X}_{tN^2} \right) \right],
\end{equation*} where $(\mathcal{X}_{tN^2})_{t \geq 0}$ is the Markov process with generator $\Delta^{2D}_N$ and with state-space {$T_N \cup \partial T_N$}. Thus, 
\begin{equation*}\begin{split}
    \sup_{t \in [0,T]} {\max_{\substack{(x,y) \in T_N}}} |\widetilde{\psi}^N_t(x,y)| &\leq \max_{\substack{(z,w) \in T_N \\ z \neq w}} |\psi^N_0(z,w)|\\& + \max_{z \in \Lambda_N} |\psi^N_0(z,z) - 1| \sup_{t \in [0,T]} {\max_{\substack{(x,y) \in T_N}}} \mathbb{P}_{(x,y)} [\mathcal{X}_{tN^2} \in \mathcal{D}_N].
\end{split}
\end{equation*}
Since the bound for {$\max_{\substack{(z,w) \in T_N}}|\psi^N_0(z,w)|$} and $\max_{z \in \Lambda_N} |\psi^N_0(z,z) - 1|$ only depends on the choice of the initial sequence of measures, i.e. {$(\mu^N)_{N \geq 1}$}, we can assume that
\begin{equation}
    {\max_{\substack{(z,w) \in T_N}}} |\psi^N_0(z,w)| \lsim \frac{1}{N} \quad \textrm{ and } \quad \max_{z \in \Lambda_N} |\psi^N_0(z,z) - 1| \leq C,
\end{equation} for some constant $C > 0$. We are now left with estimating 
\begin{equation*}
    \sup_{t \in [0,T]}{\max_{\substack{(x,y) \in T_N}}} \mathbb{P}_{(x,y)} [\mathcal{X}_{tN^2} \in \mathcal{D}_N].
\end{equation*}

We start by observing that, since $(\mathcal{X}_{tN^2})_{t \geq 0}$ is a symmetric simple random walk on $T_N \cup \partial T_N$ which has absorbing boundary, i.e. whenever $(\mathcal{X}_{tN^2})_{t \geq 0}$ hits $\partial T_N$ it can not leave it, {then, for all $(x,y) \in T_N$} and every $t \geq 0$, we have that
\begin{equation*}
    \mathbb{P}_{(x,y)} [\mathcal{X}_{tN^2} \in \mathcal{D}_N] \leq \mathbb{P}_{(x,y)} [\mathcal{X}_{tN^2} \notin \partial T_N] = \mathbb{P}_{(x,y)} [\tau_{\textrm{death}} > tN^2],
\end{equation*} where $\tau_{\textrm{death}}$ is the stopping time at which $(\mathcal{X}_{tN^2})_{t \geq 0}$ hits $\partial T_N$, i.e.
\begin{equation*}
    \tau_{\textrm{death}} := \inf \Big\{ t \geq 0 \ | \ X_{tN^2} \in \partial T_N\Big\}.
\end{equation*} We are only left to show that, for every $t \geq 0$, 
\begin{align*}
    {\max_{\substack{(x,y) \in T_N}}}\mathbb{P}_{(x,y)} [\tau_{\textrm{death}} > tN^2] \lesssim \frac{1}{\sqrt{t N^2}},
\end{align*} with a constant independent of $N$. But this is a consequence of the fact that, since $(\mathcal{X}_{tN^2})_{t \geq 0}$ is a symmetric continuous-time random walk, then it is enough to bound $${\max_{\substack{(x,y) \in T_N}}} \mathbb{P}_{(x,y)} [\tau^{1D}_{\textrm{death}} > tN^2],$$where $\tau^{1D}_{\textrm{death}}$ is the stopping time at which the projection (in one of the coordinates) of the two-dimensional random walk reaches an endpoint of the finite line $\Lambda_N$. Thus, from the classical random walk theory in one-dimension we get that whenever $t \thicksim N^2$ then $${\max_{\substack{(x,y) \in T_N}}} \mathbb{P}_{(x,y)}\Big[\tau^{1D}_{\textrm{death}} > t\Big]\lsim \frac{1}{\sqrt{t}},$$ and we are done.

We also note that the stationary solution of \eqref{equation_corr_GL}, i.e.
\begin{equation*}
\Delta^{2D}_N \psi^N_{ss}(x,y) = 2N^2 \mathbb{1}_{y=x+1} - 4 N^2 \mathbb{1}_{y=x},
\end{equation*} can be written as 
\begin{align*}
    \widetilde{\psi}^N_{ss}(x,y) = \psi^N_{ss}(x,y) - \mathbb{1}_{y=x},
\end{align*} where $\widetilde{\psi}^N_{ss}$ solves
\begin{equation*}
    \Delta^{2D}_N \widetilde{\psi}^N_{ss}(x,y) = 0,
\end{equation*} with $\widetilde{\psi}^N_{ss}(x,y)=0$ for all $(x,y) \in \partial T_N$. It is easy to see that $\widetilde{\psi}^N_{ss}(x,y) = 0$ and thus $\psi^N_{ss}(x,y) = \mathbb{1}_{y=x}$.

\subsubsection{The Brownian Energy Process (BEP)} \label{BEP}

Let $T_-, T_+ > 0$ and fix $\alpha > 0$. The Brownian Energy Process (BEP($\alpha$) or simply BEP) is a model that describes symmetric energy exchange between nearest neighboring sites. It is a Markov process $(z_t)_{t \geq 0} :=(z_t (x) \, : \, x \in \Lambda _N)_{t \geq 0}$ with 
state-space $\Omega_N=\mathbb R_+^{ \Lambda_N}$ and whose generator is given by
\begin{align} \label{generator_bep}
N^2\mcb{L}^{BEP}_N= N^2 \mcb{L}^{BEP}_- +N^2 \mcb{L}^{BEP}_{0,b} + N^2\mcb{L}^{BEP}_+,
\end{align}  which, for every $z \in \mathbb{R}^{\Lambda_N}_+$, 
\begin{align*}
\mcb{L}^{BEP}_- &= \left(T_- \alpha - \frac{1}{2} z(1) \right) \frac{\partial}{\partial z(1)} + z(1) \frac{\partial^2}{\partial z^2(1)},\\
\mcb{L}^{BEP}_+ &= \left(T_+ \alpha - \frac{1}{2} z(N-1) \right) \frac{\partial}{\partial z(N-1)} + z(N-1) \frac{\partial^2}{\partial z^2(N-1)},\\
\mcb{L}^{BEP}_{0,b} &= \sum_{x \in \Lambda_{N-1}} z(x) z(x+1) \left(\frac{\partial}{\partial z(x)} - \frac{\partial}{\partial z(x+1)}\right)^2 - \alpha (z(x) - z(x+1)) \left(\frac{\partial}{\partial z(x)} - \frac{\partial}{\partial z(x+1)} \right).
\end{align*}

The Markov process defined by the bulk generator $\mcb L^{BEP}_{0,b}$ conserves the total energy of the system
\begin{equation*}
\sum_{x \in \Lambda_N} z(x).
\end{equation*}

\paragraph{Invariant measures}
The BEP with open boundaries is irreducible and recurrent, thus it admits a unique invariant {measure}. As in \cite{CGGR18}, if $T_- = T_+ = T$, then the BEP($\alpha$) with generator $\mathcal{L}^{BEP}_N$ defined in \eqref{generator_bep} admits as unique stationary product measure the measure
\begin{align*}
\bigotimes_{x \in \Lambda_N} \textrm{Gamma}(\alpha,2T).
\end{align*} In this case, the previous measure is reversible. Out of equilibrium, i.e.  when $T_- \neq T_+$, reversibility is lost and the unique invariant measure presents non-zero long-range correlations. As in many other processes, we do not have further information about the stationary measure in this last case.

\paragraph{The correlation function}
Let {$(\mu^N)_{N \geq 1}$} be a sequence of probability measures in $\mathbb{R}_+^N$.
We define, for every $x \in \Lambda_N$ and $t \geq 0$,
$$e^N_t(x) := \mathbb{E}_{\mu^N}[z_{tN^2}(x)],$$
and we also define the two-points correlation function $\varphi^N_t$ for the BEP as, for every $x,y \in (\Lambda_N)^2$ with $x \leq y$,
\begin{equation}
\varphi^N_t(x,y) := \mathbb{E}_{\mu^N} \left[\frac{\alpha z_{tN^2}(x) \ z_{tN^2}(y)}{\alpha + \mathbb{1}_{y=x}} \right] - \mathbb{E}_{\mu^N} \left[z_{tN^2}(x) \right] \mathbb{E}_{\mu^N} \left[z_{tN^2}(y) \right].
\end{equation}

A simple computation shows that, for every $x \in \Lambda_N$ and $t \geq 0$, it holds that
\begin{align*}
\partial_t e^N_t(x) =  \Delta^{1D}_N e^N_t(x),
\end{align*} where $\Delta^{1D}_N$ coincides with the operator \eqref{laplaciannn} by taking in $c_{x,x+1}$ defined in \eqref{ratesss} the choices $c=\alpha$, $\lambda_-=T_-$ and $\lambda_+=T_+$. Also, for every $(x,y) \in T_N$,
\begin{align} \label{eq_corr_BEP}
\partial_t \varphi^N_t(x,y) = \Delta^{2D}_N \varphi^N_t(x,y) + N^2\left(e^N_t(x+1) - e^N_t(x)\right)^2 \mathbb{1}_{y=x+1},
\end{align} where $\Delta^{2D}_N$ corresponds to the operator $\Delta_N$ defined in \eqref{dual_def} by taking $d=1$, $c = \alpha$, $\lambda_-=T_-$ and $\lambda_+=T_+$. Proceeding as in Section \ref{correlation_section}, one can obtain the estimate
\begin{equation*}
\sup_{t \in [0,T]} {\max_{\substack{(x,y) \in T_N}}} |\varphi^N_t(x,y)| \lsim \frac{1}{N}.
\end{equation*}

Note that the stationary solution of \eqref{eq_corr_BEP} i.e. the solution of 
\begin{align*}
\Delta^{2D}_N \varphi^N_t(x,y) =- N^2\left(e^N_t(x+1) - e^N_t(x)\right)^2 \mathbb{1}_{y=x+1},
\end{align*} is given by
\begin{align*}
    \varphi^N_t(x,y) &= (T_+-T_{-})^2\mathcal{T}_N(x,y)\\
    &= (T_+-T_-)^2\Big(\frac{1}{N^2}\frac{(N-y)x}{N}-\frac{1}{2N^2}\textbf{1}_{x=y}\Big),
\end{align*} where $$\mathcal{T}_N(x,y) = \mathbb{E_{\mu_N}}\left[ \int_0^\infty \mathbb{1}(\mathcal{X}_{tN^2} \in \mathcal{D}^+_N)\right],$$ and it represents the occupation time of the diagonal $\mathcal{D}^+_N$ by the random walk {$(\mathcal{X}_{tN^2})_{t \geq 0}$} with generator  $\Delta^{2D}_N$.

\subsection{Interacting Piles} \label{interacting_piles_model}

Here we consider the model introduced in \cite{GiardinaFrassek} and also analyzed in \cite{Chiara_paper}. It has the particular feature of having unbound occupation variables and for each bond $\{x,x\pm 1\}$, it allows for a block of any $k \in \{1,\dots,\eta(x)\}$ particles to move from $x$ to  $x\pm 1$ with a certain rate. Next, we describe in detail the dynamics.

For real numbers $\beta_-, \beta_+ \in (0,1)$ and $\alpha \in \mathbb{N}$, we consider the continuous-time Markov process $(\eta_t)_{t \geq 0} :=(\eta_t (x) \, : \, x \in \Lambda _N)_{t \geq 0}$ with state-space $\Omega_N = (\mathbb{N} \cup \{0\})^{\Lambda_N}$ and  whose  infinitesimal generator $\mcb{L}_N^{IP}$  is given by
\begin{equation*}
    N^2\mcb{L}_N^{IP} := N^2\mcb{L}^{IP}_- + N^2\mcb{L}^{IP}_{0,b}+ N^2\mcb{L}^{IP}_+,
\end{equation*} 
whose action on local functions $f: \Omega_N \to \mathbb{R}$ can be written as
\begin{align*}
    (\mcb{L}^{IP}_- f)(\eta) &:= \sum_{k=1}^{\eta(1)} h_\alpha(k,\eta(1))[f(\eta - k \delta_1) - f(\eta)] + \sum_{k=1}^\infty \frac{\beta_-^k}{k} [f(\eta + k \delta_1) - f(\eta)], \\
    (\mcb{L}^{IP}_{0,b} f)(\eta) &:= \sum_{x \in \Lambda_{N-1}} \sum_{k=1}^{\eta(x)} h_\alpha(k,\eta(x))[f(\eta - k \delta_x + k \delta_{x+1}) - f(\eta)] \\
    &+ \sum_{x \in \Lambda_{N-1}} \sum_{k=1}^{\eta(x+1)} h_\alpha(k,\eta(x+1))[f(\eta - k\delta_{x+1} + k\delta_{x}) - f(\eta)], \\
    (\mcb{L}^{IP}_+ f)(\eta) &:= \sum_{k=1}^{\eta(N-1)} h_\alpha(k,\eta(N-1))[f(\eta - k \delta_{N-1}) - f(\eta)] + \sum_{k=1}^\infty \frac{\beta_+^k}{k} [f(\eta + k \delta_{N-1}) - f(\eta)], 
\end{align*} where the function $h_\alpha: {\mathbb{N} \times \mathbb{N}} \to \mathbb{R}$ is given by
\begin{equation}
    {h_\alpha(j,m) = \frac{1}{j} \frac{\Gamma(m+1) \Gamma(m-j+ \alpha)}{\Gamma(m-j+1) \Gamma(m+\alpha)} \mathbb{1}_{1 \leq j \leq m}},
\end{equation} and  $\Gamma(\cdot)$ represents the gamma function.

The Markov process defined by the bulk generator $\mcb L^{IP}_{0,b}$ conserves the total number of particles
\begin{equation*}
\sum_{x \in \Lambda_N} \eta (x).
\end{equation*}

\subsubsection{Invariant measures}

Because this process is irreducible and recurrent, it has a unique invariant measure. 

\begin{itemize}
    \item If $\beta_- = \beta_+ = \beta$, i.e. in the equilibrium case, the unique invariant measure $\mu_{ss}$ is given by a product of Negative Binomial distributions with parameters $\alpha > 0$ and $0 < \beta < 1$. Namely
    \begin{align*}
        \mu_{ss}(\eta) = \prod_{x \in \Lambda_N} \frac{\beta^{\eta(x)}}{\eta(x)!} \frac{\Gamma(\eta(x) + \alpha)}{\Gamma(\alpha)}(1- \beta)^\alpha.
    \end{align*}

    \item If $\beta_- \neq \beta_+$, i.e. in the non-equilibrium case, reversibility is lost and one can check that a product ansatz for the stationary measure does not work. Indeed, the non-equilibrium steady state has long-range correlations as it is shown in \cite{GiardinaFrassek}.
\end{itemize}

\subsubsection{The correlation function}

Let {$(\mu^N)_{N \geq 1}$} be a sequence of probability measures in {$(\mathbb{N} \cup \{0\})^{\Lambda_N}$}. If we define the discrete density for the model defined in Section \ref{interacting_piles_model} as
$$\rho^N_t(x) := \mathbb{E}_{\mu^N}[\eta_{tN^2}(x)],$$ and also define its two-points  correlation function $\varphi^N_t$ as
\begin{equation}
\varphi^N_t(x,y) := \mathbb{E}_{\mu^N}\left[\frac{\alpha \eta_{tN^2}(x)[\eta_{tN^2}(y)  {-} \mathbb{1}_{y=x}]}{\alpha + \mathbb{1}_{y=x}} \right] - \mathbb{E}_{\mu^N}[\eta_{tN^2}(x)] \mathbb{E}_{\mu^N}[\eta_{tN^2}(y)],
\end{equation} 
for every {$x \in \Lambda_N$}, a simple computation shows that
\begin{align} \label{density_equation}
\partial_t \rho^N_t(x)  = {\Delta}^{1D}_N \rho^N_t(x),
\end{align} where ${\Delta}^{1D}_N$ coincides with the operator \eqref{laplaciannn} by taking the choice for  {$c_{x,y}$ where $|x-y|=1$ with $c=1/\alpha$, $\lambda_-=1$ and $\lambda_+=1$. Here $\rho^N_t(0) := \frac{\beta_-}{1-\beta_-}$ and $\rho^N_t(N) := \frac{\beta_+}{1-\beta_+}$}. Moreover, for every $x,y \in (\Lambda_N)^2$ with $x \leq y$,
\begin{align} \label{equation_corr_pile}
\partial_t \varphi^N_t(x,y) =  {N^2 \mathcal{L}_N\varphi^N_t(x,y) + [N^2(\rho^N_t(x+1) - \rho^N_t(x))^2 + N^2(\rho^N_t(x-1) - \rho^N_t(x))^2] \mathbb{1}_{y = x}},
\end{align} where in \eqref{equation_corr_pile}  {$\mathcal{L}_N$ corresponds to the operator 
 given, for every function $f: (\Lambda_N \cup \{0,N\})^2 \to \mathbb{R}$ with $f(x,y) = 0$ for every $(x,y) \in \partial \Lambda_N$, by}
 {\begin{align*}
\mathcal{L}_N f(x,y) = \begin{cases}
\frac{1}{\alpha}[f(x-1,y)+f(x+1,y)+f(x,y+1)+f(x,y-1)-4f(x,y)], \\
\qquad\textrm{ if } |y-x| \geq 1,\\
\frac{1}{\alpha(\alpha+1)}[2f(x-1,x)+2f(x,x+1) +f(x-1,x-1)+f(x+1,x+1)-6f(x,x)],\\
\qquad \textrm{ if } y=x.
\end{cases}
\end{align*}} Proceeding  {with a similar strategy as the one } in Section \ref{correlation_section}, one can obtain that
\begin{equation*}
\sup_{t \in [0,T]} {\max_{\substack{(x,y) \in T_N}}} |\varphi^N_t(x,y)| \lsim \frac{1}{N}.
\end{equation*}

 {We also note that the stationary solution of \eqref{density_equation}, i.e.
\begin{equation*}
    {\Delta}^{1D}_N \rho^N_{ss}(x) = 0,
\end{equation*} is given by
\begin{equation}
    \rho^N_{ss}(x) = \left(\frac{\beta_+}{1-\beta_+} - \frac{\beta_-}{1-\beta_-}\right)\frac{x}{N} + \frac{\beta_-}{1-\beta_-}.
\end{equation} Moreover, the stationary solution} of \eqref{equation_corr_pile} i.e. the solution of 
\begin{align*}
 {N^2\mathcal{L}_N \varphi^N_{ss}(x,y) = - N^2[(\rho^N_{ss}(x+1) - \rho^N_{ss}(x))^2 + (\rho^N_{ss}(x-1) - \rho^N_{ss}(x))^2] \mathbb{1}_{y = x},}
\end{align*} is given by 
 {\begin{equation*}
\begin{split}
\varphi_{ss}^N(x,y)&=\frac{\alpha(\alpha+1)}{N^2(N+1)}\left(\frac{\beta_+}{1-\beta^+}-\frac{\beta_-}{1-\beta_-}\right)^2(N-y)x,
\end{split}
\end{equation*}} and thus
\begin{equation*}
{\max_{\substack{(x,y) \in T_N}}} |\varphi^N_{ss}(x,y)| \lsim \frac{1}{N}.
\end{equation*}

\appendix

\section{Remark on stochastic duality} 	\label{remark on duality}

\subsection{On the models of interacting particles}

Following Section 4.1 of \cite{CGGR18}, it is well known that both (SEP), (SIP), and (IRW) with an open boundary have a dual process with a pure absorbing boundary through a duality function, which is a function of product type. More precisely, setting $M \in \{(SEP), (SIP), (IRW)\}$ for the respective choices of $c$ and $d$ mentioned in Section \ref{sec_dynamics}, the model $M$ with an open boundary has the model $M$ with absorbing boundary (and state-space $\widehat{\Omega}^{M}_N$) as dual process, with a duality function $D^{M}: \Omega^{M}_N \times \widehat{\Omega}^{M}_N \to \mathbb{R}$ given, for every $\eta \in \Omega^M_N$ and $\widehat{\eta} \in \widehat{\Omega}^M_N$, by
\begin{align}
    D^{M}(\eta,\widehat{\eta}) = \left(\frac{\varrho_-}{c}\right)^{\widehat{\eta}(0)} d^M(\eta(x),\widehat{\eta}(x)) \left(\frac{\varrho_+}{c}\right)^{\widehat{\eta}(N)},
\end{align}
where
\begin{align*}
    \Omega^M_N = \begin{cases} \{0,\dots,c\}^{\Lambda_N}, \textrm{ if } M = (SEP),\\
    (\mathbb{N} \cup \{0\})^{\Lambda_N}, \textrm{ if } M = (SIP) \textrm{ or } M= (IRW),
    \end{cases}
\end{align*}
\begin{align*}
\widehat{\Omega}^M_N = \begin{cases} (\mathbb{N}\cup\{0\}) \times \{0,\dots,c\}^{\Lambda_N} \times (\mathbb{N}\cup\{0\}) , \textrm{ if } M = (SEP),\\
(\mathbb{N} \cup \{0\})^{\Lambda_N \cup \{0,N\}}, \textrm{ if } M = (SIP) \textrm{ or } M= (IRW),
\end{cases},
\end{align*} and
\begin{align*}
d^M(\eta(x),\widehat{\eta}(x)) = \begin{cases}
    \prod_{x \in \Lambda_N} \frac{\eta(x)!}{(\eta(x) - \widehat{\eta}(x))!}\frac{(c - \hat{\eta}(x))!}{c!}, \textrm{ if } M = (SEP),\\
    \vspace{0,01cm}\\
    \prod_{x \in \Lambda_N} \frac{\eta(x)!}{(\eta(x) - \hat{\eta}(x))!}\frac{\Gamma(c)}{\Gamma(c+\hat{\eta}(x))}, \textrm{ if } M = (SIP),\\
    \vspace{0,01cm}\\
    \prod_{x \in \Lambda_N} \frac{\eta(x)!}{(\eta(x) - \hat{\eta}(x))!}, \textrm{ if } M = (IRW).
    \end{cases}
\end{align*} Here $M$ with absorbing boundary is the process with Markov generator 
\begin{equation} \label{dual_sep_sip_irw}
\widehat{\mcb{L}}^M_{N} := N^2(\mcb{L}_{0,N} + \widehat{\mcb{L}}^M_- + \widehat{\mcb{L}}^M_+),
\end{equation}where, for $j \in \{ -, +\}$, we have that, for every $\widehat{\eta} \in \widehat{\Omega}_N$ and all $f: \widehat{\Omega}_N \to \mathbb{R}$,
\begin{align} \label{dual_def_bound}
   \widehat{\mcb{L}}^M_j f(\widehat{\eta}) := c \lambda_j \widehat{\eta}(x^j) \big\{ f(\widehat{\eta} - \delta_{x^j} + \delta_{y^j}) -f(\widehat{\eta}) \big\},
\end{align} {where $x^- = 1$ and $x^+ = N-1$ and } with $y^j = 0$ if $x^j = 1$ and $y^j = N$ if $x^j = N-1$.

Using the expression of $D^M$ given above, we obtain that
\begin{equation}
\varphi^N_t(x,y) = c^2 \left[ \mathbb{E}_{\mu^N}[ D^M (\cdot,\delta_x + \delta_y)] - \mathbb{E}_{\mu^N}[D^M (\cdot,\delta_x)] \mathbb{E}_{\mu^N}[D^M (\cdot,\delta_y)] \right].
\end{equation} Thus, for every $(x,y) \in T_N$,
\begin{align*}
\partial_t \varphi^N_t(x,y) = \widehat{\mcb{L}}^M_N \varphi^N_t(x,y) + d N^2\left(\rho^N_t(x+1) - \rho^N_t(x)\right)^2 \mathbb{1}_{y=x+1},
\end{align*}  where   $\rho^N_t(x)=\mathbb E_{\mu^N}[\eta_{tN^2}(x)]$. From \eqref{eq_varphi} we conclude that $\Delta^{2D}_N = \widehat{\mcb{L}}^M_N$.

\subsection{Other models}
\subsubsection{BEP}

It is also known that (BEP) with open boundary and (SIP) with absorbing boundary (taking as boundary parameters in \eqref{dual_def_bound} $\lambda_-=\lambda_+=\frac{1}{2\alpha}$) are dual processes with duality function $D^{BEP}$ given by
\begin{align*}
D^{BEP} (z, \xi) := (2 T_-)^{\xi(0)} \left[\prod_{x \in \Lambda_N} [z(x)]^{\xi(x)} \frac{\Gamma (\alpha)}{\Gamma(\alpha + \xi(x))} \right] (2 T_+)^{\xi(N)},
\end{align*} as is defined in \cite{CGGR18} equation (4.9),  where $T_-$ and $T_+$ are the constants introduced in the beginning of Section \ref{BEP}; and in the case without boundary it is defined below equation (16) of \cite{Redig_Sau2018}. Here the generator of the dual, i.e. (SIP) with absorbing boundary, that we denote by $\mcb{L}^{dual}_{SIP}$, is given by {\eqref{dual_sep_sip_irw} where in \eqref{dual_def_bound} we take $\lambda_-=\lambda_+=\frac{1}{2\alpha}$}. Using the expression of $D^{BEP}$ given above, we obtain that
\begin{equation}
\varphi^N_t(x,y) = {\alpha^2 \left[\mathbb{E}_{\mu^{BEP}_N}[ D^{BEP} (z,\delta_x + \delta_y)] - \mathbb{E}_{\mu^{BEP}_N}[D^{BEP} (z,\delta_x)] \mathbb{E}_{\mu^{BEP}_N}[D^{BEP} (z,\delta_y)] \right].}
\end{equation}  

Also, for every $(x,y) \in T_N$ with $x \leq y$,
\begin{align*}
\partial_t \varphi^N_t(x,y) = \mcb{L}^{dual}_{SIP} \varphi^N_t(x,y) + N^2\left(e^N_t(x+1) - e^N_t(x)\right)^2 \mathbb{1}_{y=x+1}.
\end{align*} 
Above $e^N_t(x) := \mathbb{E}_{\mu^N}[z_{tN^2}(x)]$. From \eqref{eq_corr_BEP}, we have that the operator  $\Delta^{2D}_N$ which corresponds to the operator  defined in \eqref{dual_def} taking $d=1$, $c = \alpha$, $\lambda_-=T_-$ and $\lambda_+=T_+$, is such that  ${\Delta^{2D}_N} = \mcb {L}^{dual}_{SIP}$.

\subsubsection{Interacting Piles}

Regarding the model introduced in Section \ref{interacting_piles_model}, it was shown in \cite{GiardinaFrassek} that this model is dual to its version with only absorbing boundaries  - see Definition 2.3 of \cite{GiardinaFrassek}. In the dual process, when particles are in the bulk they move as in the original system; but when they reach the boundary point $0$ (resp. $N$), they are absorbed at a rate, for each $k \in \{1,\dots,\eta(1)\}$ (resp., $k \in \{1,\dots,\eta(N-1)\}$), $h_\alpha(k,\eta(1))$ (resp., $h_\alpha(k,\eta(N-1))$). {More precisely, denoting} $\widehat{\Omega}_N := (\mathbb{N} \cup \{0\})^{\Lambda_N \cup \{0,N\}}$, through the duality function $D: {\Omega_N} \times \widehat{\Omega}_N \to \mathbb{N}$, which is given, for every $(\eta,\xi) \in {\Omega_N} \times \hat{\Omega}_N$,
\begin{equation*}
    D(\eta, \xi) = (\rho_-)^{\xi(0)} \left[\prod_{x \in \Lambda_N} \frac{\eta(x)!}{(\eta(x) - \xi(x))!} \frac{\Gamma(\alpha)}{\Gamma(\alpha+\xi(x))} \right](\rho_+)^{\xi(N)},
\end{equation*} where here $\rho_\pm= \frac{\beta_\pm}{1-\beta_\pm}$, the model with generator $\mcb{L}^{IP}$ is dual to the continuous-time Markov process with state-space $\widehat{\Omega}_N$ and Markov generator $\mcb{L}^{dual}_{IP}$ given by
\begin{align*}
    \mcb{L}^{dual}_{IP} = N^2(\mcb{L}^{IP,dual}_- + \mcb{L}^{IP}_{0,b} + \mcb{L}^{IP,dual}_+),
\end{align*}
where, for every $f :\widehat{\Omega}_N\to\mathbb R$ and $\xi \in \widehat{\Omega}_N$, 
\begin{align*}
    (\mcb{L}^{IP,dual}_- f)(\xi) = \sum_{k=1}^{\xi(1)} h_\alpha(k,\xi(1))[f(\xi - k \delta_1 + k \delta_0) - f(\xi)],
\end{align*} and, analogously,
\begin{align*}
    (\mcb{L}^{IP,dual}_+ f)(\xi) = \sum_{k=1}^{\xi(N-1)} h_\alpha(k,\xi(N-1))[f(\xi - k \delta_{N-1} + k \delta_N) - f(\xi)].
\end{align*} From the definition of $D$, it is easy to verify that
\begin{equation*}
    \rho^N_t(x) = \alpha \mathbb{E}_{\mu^N}[D(\eta_{tN^2}, \delta_x)],
\end{equation*} and 
\begin{equation*}
    \varphi^N_t(x,y) = {\alpha^2 \left[ \mathbb{E}_{\mu^N}[D(\eta_{tN^2}, \delta_x+\delta_y)] - \mathbb{E}_{\mu^N}[D(\eta_{tN^2}, \delta_x)] \mathbb{E}_{\mu^N}[D(\eta_{tN^2}, \delta_y)] \right]}.
\end{equation*} Moreover
\begin{equation*}
    \partial_t \rho_t^N(x) = \mcb{L}^{dual}_{IP} \rho^N_t(x),
\end{equation*} and
\begin{equation*}
    \partial_t \varphi^N_t(x,y) = \mcb{L}^{dual}_{IP} \varphi^N_t(x,y) + N^2\left(\rho^N_t(x+1) - \rho^N_t(x)\right)^2 \mathbb{1}_{y=x+1}.
\end{equation*} Thus the operator {$N^2\mathcal{L}$} given in \eqref{equation_corr_pile} corresponds to $\mcb{L}^{dual}_{IP}$.

\begin{acknowledgements}
B.Salvador thanks FCT/Portugal for the financial support through the PhD scholarship with reference 2022.13270.BD. P. Gonçalves thanks Funda\c c\~ao para a Ci\^encia e Tecnologia FCT/Portugal for financial support through the projects UIDB/04459/2020, UIDP/04459/2020 and the project SAUL funded by FCT-ERC.
\end{acknowledgements}

\end{document}